\documentclass[12pt]{article}

\usepackage[english, main=english]{babel}
\usepackage{epsfig}
\usepackage{dsfont}
\usepackage{subfig}
\usepackage{amsmath,amssymb,amsthm,thmtools}
\usepackage{tikz}
\usepackage{xparse}
\usepackage{ifthen}
\usepackage{textcomp}
\usepackage{mathtools}
\usepackage{mathrsfs}
\usepackage{enumitem}
\usepackage{comment}
\usepackage[T1]{fontenc}
\usepackage[utf8]{inputenc}
\usepackage{graphicx}
\usepackage{pdfpages}
\usepackage{float}
\RequirePackage{times} 
\usepackage[colorlinks,linkcolor=blue,urlcolor=black]{hyperref}
\usepackage[dvipsnames]{xcolor}

\newcommand{\N}{\mathbb{N}} 
\newcommand{\R}{\mathbb{R}} 
\newcommand{\C}{\mathbb{C}} 
\renewcommand{\P}{\mathbb{P}} 

\renewcommand{\d}{\mathrm{d}}

\theoremstyle{definition}
\newtheorem{definition}{Definition}[section]


\theoremstyle{plain}

\newtheorem{Lemma}[definition]{Lemma}
\newtheorem{Corollary}[definition]{Corollary}
\newtheorem{Proposition}[definition]{Proposition}
\newtheorem{Theorem}[definition]{Theorem}

\theoremstyle{remark}
\newtheorem{Remark}{Remark}

\title{MA(1) processes with uniform innovations conditioned to stay positive in the expanding regime}
\author{Frank Aurzada and  Virginia Worf}

\begin{document}

\maketitle

\begin{abstract}
We study MA(1) processes with uniform innovations conditioned to remain positive for a long time when $|\theta|>1$, where $\theta$ denotes the coupling parameter. We show that, for a fixed starting point, the finite-dimensional distributions of the conditioned process converge and identify the limiting process as a Doob $h$-transform with explicit function $h$. We also obtain precise asymptotics for the corresponding persistence probabilities. The limiting process exhibits qualitatively different behaviour depending on $\theta$ and on the support of the innovation distribution. We also determine stationary distributions of the limiting conditioned processes for all $\theta\ne 1$. 
\end{abstract}
\section{Introduction and main results}
Conditioning stochastic processes to be positive over a long time is a classical problem in probability going back to the seminal work of Doob \cite{doob}. It is part of the more general area of stochastic processes under constraints and closely related to the questions of persistence and quasi-stationarity. Apart from finding the probability of the conditioning event, one would like to understand how the conditioning affects the dynamics of the process. Moving-average processes with uniform innovations provide a simple model for which both the transition probabilities and the stationary distribution of the limiting Markov chain are non-trivial but can be determined fully explicitly. This makes it possible to see how the coupling parameter and the support of the innovation distribution shape the long-term behavior
of the conditioned process. 

We focus on MA(1) processes in the expanding regime $|\theta|>1$,
where $\theta$ is the coupling parameter. For a fixed starting point, we obtain sharp asymptotics for the persistence
probabilities and show that the conditioned finite-dimensional
distributions converge to those of a Doob $h$-transform for an explicit function $h$.
We also determine the stationary distributions of the limiting conditioned chains, including those in the non-expanding regime $\theta\in[-1,1)$.

The non-expanding case was studied in \cite{AW}.
The two cases are related via a time-reversal argument; however, the results for the non-expanding regime do not directly determine the limiting conditioned process in the expanding case.  In particular, the function $h$ defining the Doob transform cannot generally be recovered from the corresponding function in the dual non-expanding regime. Its explicit determination, together with the proof of the fixed-starting-point persistence asymptotics and the explicit form of the stationary distribution, form the main contributions of the present paper. 

\begin{figure}[H]
\begin{center}
\begin{tikzpicture}[scale=0.8]
\filldraw[fill=blue!50!white, draw=blue!50!white] (-4,4) -- (0,0) -- (0,4) -- cycle; 
\filldraw[fill=blue!50!white, draw=blue!50!white] (0,4) -- (0,0) -- (4,0) -- (4,4) -- cycle; 


\filldraw[fill=blue!30!white, draw=blue!30!white] (4,0) rectangle (8,4); 

\filldraw[fill=blue!50!white, draw=blue!50!white] (0,4) rectangle (4,8); 


\filldraw[fill=blue!30!white, draw=blue!30!white] (4,4) rectangle (8,8);

\fill[blue!50!white]  plot[domain=-4:-2, samples=100] (\x, {-16/\x})   -- (-2,4) -- (-4,4) -- cycle;
\fill[blue!50!white]
  (0,4)  -- (0,8)  -- (-2.1,8)  -- (-2.1,4)  -- cycle;


\fill[green!60!white]
  plot[domain=-4:-2, samples=100] (\x, {-16/\x})
  -- (-2, 8) -- (-4, 8) -- cycle;

\fill[green!30!white]
  (-8,4)
  -- (-8,2)
  plot[domain=-8:-4, samples=100] (\x, {-16/\x})
  -- (-4,4)
  -- (-8,4)
  -- cycle;

\filldraw[fill=yellow!30!white, draw=yellow!30!white] (-4,4) -- (-5,3.2) -- (-6,2.67) -- (-7,2.29) -- (-8,2) -- (-8,0) -- (-4,0) -- cycle;

\filldraw[fill=yellow!60!white, draw=yellow!60!white] (-4,4) -- (0,0) -- (-4,0) -- cycle;


\filldraw[fill=orange!80!white, draw=orange!30!white] (0,0) -- (4,-4) -- (4,0) -- cycle;



\fill[gray!30!white]
  plot[domain=4:8, samples=100] (\x, {-16/\x})
  -- (8,0) -- (4,0) -- cycle;
\draw[-,dashed] (4,8) -- (4,-4);
\draw[-,dashed] (-4,8) -- (-4,0);


\node [black] at (-4,-2) {$p_n^{a,y}(\theta)\equiv 1$};
\node [black] at (1.5,-3.3) {$p_n^{a,y}(\theta)\equiv 1$};
\node [black] at (7,-3.3) {$p_n^{a,y}(\theta)\equiv 0,n\ge 1$};
\node[black, align=center] at (6,-1.2){$p_n^{a,y}(\theta)\equiv0$\\for large $n$};

\filldraw[fill=green!30!white, draw=green!40!white] (-8,8) -- (-8,4) -- (-4,4) -- cycle;
\filldraw[fill=green!30!white, draw=green!40!white] (-8,8) -- (-4,8) -- (-4,4) -- cycle;

\draw[-,dashed] (8,4) -- (-8,4);

\draw[->,very thick] (-9,0) -- (9,0);
\draw[-,very thick] (0,-4) -- (0,0);
\draw[->, RoyalBlue, very thick] (0,0) -- (0,9);
\node [black] at (9.5,0) {\large $\theta$};
\node [black] at (4,-0.5) {\large $1$};
\node [black] at (-4,-0.5) {\large $-1$};
\node [black] at (0,9.5) {\large $a$};
\node [black] at (0.3,4.4) {\large $1$};
\node [black] at (0.5,-0.5) {\large $0$};
\node [black] at (-0.7,-3.7) {\large $-1$};
\draw[-,thick] (-9,-4) -- (9,-4);
\node[black] at (6,6) {\textbf{(B)}};
\node[black] at (-6,5.2) {\textbf{(G)}};
\node[black] at (-6,1) {\textbf{(Y)}};
\node[black] at (-6,-2) {\textbf{(W)}};
 \filldraw[fill=orange!30!white, draw=orange!30!white] (4,-0.04) -- (8,-0.04) -- (8,-0.04) -- (4,-0.04) -- cycle;
\end{tikzpicture}
\end{center}
\caption{Phase diagram. The light regions with $|\theta|>1$ are treated here, while the dark regions with $|\theta|<1$ were treated in \cite{AW}. The light blue region for $\theta<-1$ in
\cite{AurzadaRaschel} is included in our light green
region, since the same effective state space and
construction of $h$ apply on both sides of $a=-\theta$. }
\label{fig:map}
\end{figure}
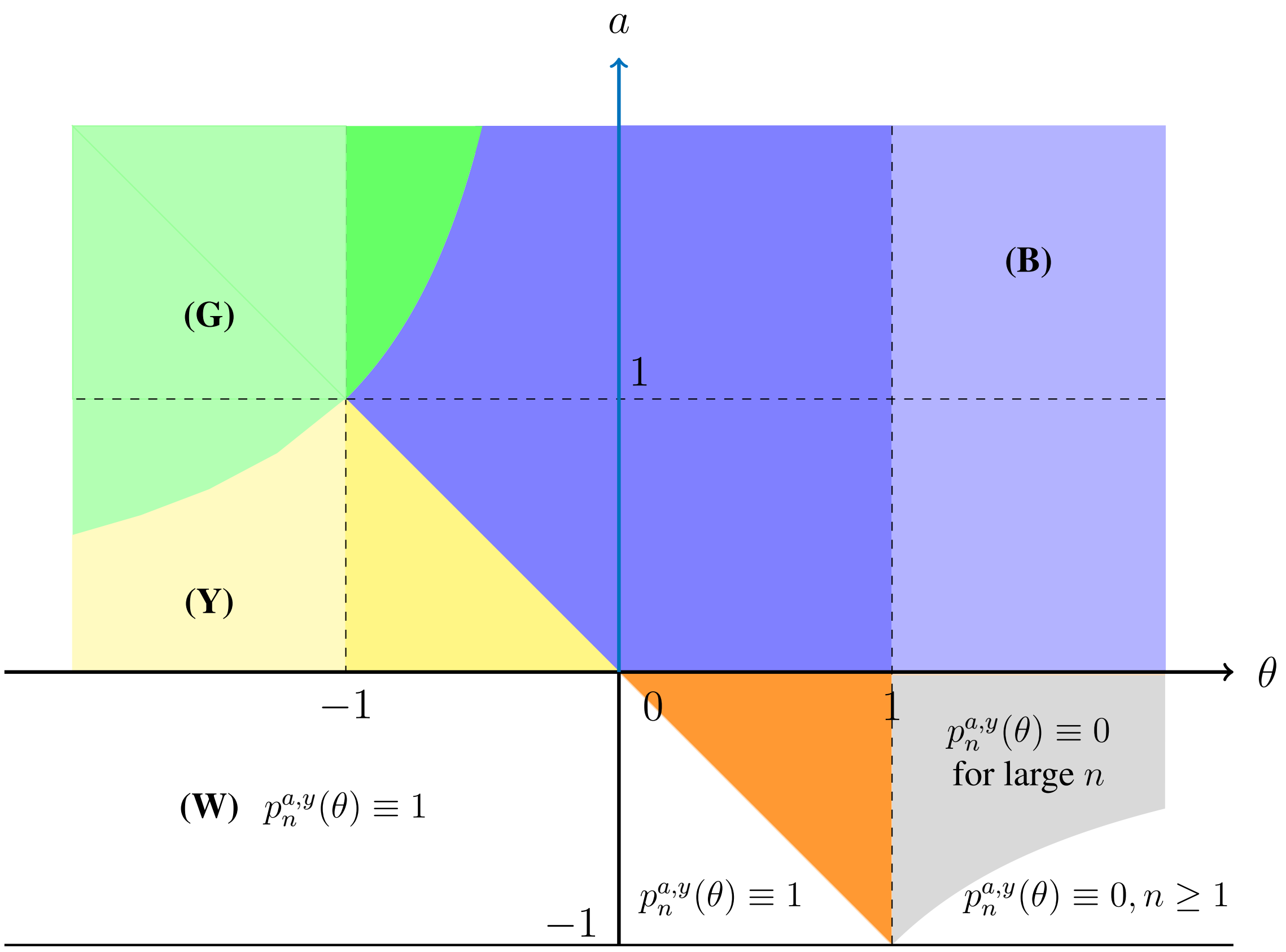

Our results show that the limiting Markov chain exhibits different structures depending on the parameters, cf.\ Figure~\ref{fig:map}. The function $h$ can be written down explicitly in all the different regimes. Interestingly, for $\theta<-1$, the function $h$ is infinitely differentiable at zero, although its Taylor series has radius of convergence zero.

More concretely, we consider the MA(1)-process  $Z_i\coloneqq X_{i+1}-\theta X_i$, $i\in\N$, $Z_0=X_1$, where $(X_i)_{i\in\N}$ are i.i.d.\ random variables, uniformly distributed on $[-a,1]$ for some $a>-1$ and $\theta\in\R$ is a coupling parameter. We normalize the positive right endpoint of the uniform distribution
to $1$. This does not restrict the generality of our results. The process $Z$ can be described using the Markov chain defined by $M_0:=(X_1,0)$ and  $M_i\coloneqq(X_{i+1},X_i)$ for $i\geq 1$. With $$S\coloneqq \{(x_2,x_1)\in [-a,1]^2| x_2>\theta x_1\}$$ and $\Omega_n\coloneqq \cap_{i=1}^n\{M_i\in S\}$, the persistence probability of the process $Z$ is given by
    \begin{align*}
    p_n^a(\theta)\coloneqq \P(\min_{1\leq i\leq n} Z_i\geq 0)=\P(X_2\geq\theta X_1,X_3\geq\theta X_2,\dots,X_{n+1}\geq\theta X_n)=\P(\Omega_n),~~~~n\ge 1,
\end{align*} 
with $p_0^a(\theta)\coloneqq 1$. Let $\P_y(~\cdot~)\coloneqq\P(~\cdot~|X_1=y)$. Then the limiting finite dimensional distributions of $(Z_i)_{i\in\N}$ conditioned to stay positive and started at $Z_0=y$ are given by
\begin{align}\label{def_lim_process}
    \lim_{n\to\infty}\P_y((M_k)_{1\leq k\leq m}\in\cdot~|\Omega_n)=\P_y((M'_k)_{1\leq k\leq m}\in\cdot~),
\end{align}
provided that a process $M'$ on the right-hand side exists. 
We prove that a Markov chain $M'$ exists and give an explicit formula for its transition kernel. Before stating the main theorem, we introduce some notation.

The transition kernel of the Markov chain $M$ for $(x_1,x_1'),(x_2,x_2')\in[-a,1]^2$ is given by
\begin{align*}
     \P_{(x_1,x_1')}(M_1\in (\d x_2,\d x_2'))=\delta_{x_1}(\d x_2')  \P(X_1\in \d x_2)=\delta_{x_1}(\d x_2')  \frac{1}{a+1}\d x_2.
\end{align*} 
Since $x_2'=x_1$ and $x_1'$ is irrelevant, it is enough to consider $x_1,x_2\in[-a,1]$. We therefore use the reduced transition kernel
 \begin{align*}
       p(x_1,\d x_2)\coloneqq \frac{1}{a+1}\d x_2.
\end{align*}
In the following, we will use the notation 
    \begin{align*}
    p_n^{a,y}(\theta)\coloneqq\P_{y}(\Omega_n)=\P(X_2\geq\theta y,X_3\geq\theta X_2,\dots,X_{n+1}\geq\theta X_n),~~~~n\ge 1,
    \end{align*}
 with $p_0^{a,y}(\theta)\coloneqq1$  for the persistence probability of the Markov chain $(M_i)_{i\in\N_0}$ staying in $S$ and started at $M_0=(y,0)$ with $y\in[-a,1]$. For $\theta\neq0$, we write $\bar\theta\coloneqq\theta^{-1}$ and for $x\in\R$, write $x_+\coloneqq\max\{x,0\}$.

To describe our main results, we need the  deformed exponential function (see \cite{Morris} for background): For $r,z\in\C$ with $|r|\leq1$, with the convention $0^0=1$, set
\begin{align*}
E(r,z)\coloneqq \sum_{n=0}^\infty \frac{r^{n(n-1)/2}}{n!}z^n.
\end{align*}
We can now state the main theorem. Depending on the parameters, the positivity constraint confines
the limiting conditioned chain to a proper subset of $[-a,1]$.
We denote this effective state space by $\Omega$ and specify it
in each case. We distinguish four parameter regimes for $a$ and the coupling parameter $\theta$, corresponding to one of the regions light blue, light green, light yellow and white presented in Figure~\ref{fig:map} below. The darker regions with $|\theta|<1$ were studied in \cite{AW} and will be referred to in Section~\ref{sec_prelim}. The remaining expanding cases are discussed in the remarks below. 
\begin{Theorem}\label{maintheorem_ge1}
    Let $|\theta|>1$. Assume that $a>-1$ if $\theta<0$ and
    $a>0$ if $\theta>0$ and consider $\Omega\subseteq[-a,1]$ given in the different regions below and $y\in \Omega$. Then, the limiting process $M'$ in \eqref{def_lim_process} exists. Moreover, it is a time-homogeneous Markov chain with reduced transition kernel 
\begin{align*}
     p'(x_1,\d x_2)&\coloneqq\mathds{1}_S(x_2,x_1)\frac{h(x_2)}{h(x_1)}\frac{1}{\lambda}\frac{1}{a+1}\d x_2,~~~~~~~~x_1,x_2\in \Omega,
\end{align*}
 where $h:\R\to\R$ is continuous and non-negative, strictly positive on
    $\Omega$, equal to zero on $[-a,1]\setminus\Omega$ and satisfies
 \begin{align}\label{eq_eigenfunction}
        \lambda h(x)=\int_{\{y>\theta x\}\cap \Omega}
        h(y)\frac{1}{1+a}\d y,
        \qquad x\in\R.
    \end{align}
In cases (B), (G), and (Y), there is a constant $c(\theta,\lambda,a)>0$, independent of $y$, such that
\begin{align*}
    p_n^{a,y}(\theta)
    \sim c(\theta,\lambda,a)h(y)\lambda^{n+1},
    \qquad y\in\Omega.
\end{align*}
    The values of $\Omega$, $\lambda$ and $h$ in the four regions are as follows.
\begin{enumerate}[label={\rm (\Alph{*})},ref={\rm (\Alph{*})}]
\setcounter{enumi}{1}
        \item Let $\theta>1$ and $a>0$. Then, $\Omega=[-a,0)$, $\lambda$ is the inverse of the smallest positive zero of $ z\mapsto
            E\left(\bar\theta,\frac{-az}{1+a}\right)$ and, $h(x)=0$ for $x\ge 0$, while for $x<0$
        \begin{align*}
            h(x)\coloneqq1+\sum_{m=1}^{\infty}\frac{\left(\frac{-\theta}{\lambda(1+a)}\right)^m\theta^{m(m-1)/2} }{m!}\left(a\bar\theta^m+x\right)_+^m.
        \end{align*}
\setcounter{enumi}{6}
        \item Let $\theta<-1$ and $a\ge-\bar\theta$. Then, $\Omega=(\bar\theta,1]$, $\lambda$ is the inverse of the smallest positive zero of $z\mapsto E\left(\bar\theta,\frac{-z}{1+a}\right)$ and the restriction of $h$ to $[-a,1]$ is the continuous function which vanishes on $[-a,\bar\theta]$, is normalized by $h(x)=1$ on $\bar\theta^2\leq x\leq1$, and is determined on $(\bar\theta,\bar\theta^2)$ by $h(\bar\theta^2)=1$ and 
        \begin{align*}
            h'(x)=\frac{-\theta}{\lambda(1+a)} h(\theta x).
        \end{align*}
        Equivalently, $h$ is given by the explicit sum formulas in Lemma~\ref{method_of_steps_theta_smaller_minus_one_operator}.
 \setcounter{enumi}{24}      
        \item Let $\theta<-1$ and $0<a<-\bar\theta$. Then, $\Omega= [-a,1]$, $\lambda$ is the inverse of the smallest positive zero of $z\mapsto E\left(\bar\theta,
    \frac{za}{\bar\theta(1+a)}\right)
    -z\frac{\bar\theta+a}{\bar\theta(1+a)}
    E\left(\bar\theta,\frac{za}{1+a}\right)$ and the restriction of $h$ to $[-a,1]$ is the continuous function normalized by $h(x)=1$ on $-a\bar\theta\leq x\leq1$, with $h(-a)=\frac{1+\theta a}{\lambda(1+a)}$,
        and determined on $(-a,-a\bar\theta)$  by $h(-a\bar\theta)=1$ and 
        \begin{align*}
            h'(x)=\frac{-\theta}{\lambda(1+a)} h(\theta x).
        \end{align*}
        Equivalently, $h$ is given by the explicit sum formulas in Lemma~\ref{method_of_steps_theta_smaller_minus_one_operator}.
    \setcounter{enumi}{22}
        \item Let $\theta<-1$ and $-1<a\leq0$. Then, $\Omega=[-a,1]$, $\lambda=1$ and $h(x)=1$ for $x\in[-a,1]$.
    \end{enumerate}
\end{Theorem}
\begin{Remark}\label{Case_B_y_pos}
In Theorem~\ref{maintheorem_ge1} (B), where $\theta>1$ and $a>0$, the restriction to $\Omega=[-a,0)$ is necessary. Indeed, for $y>0$ we have $\Omega_n=\emptyset$ as soon as $\theta^ny\geq1$. Thus, a positive starting point leads to the event conditioned on to be empty for large enough $n$. Similarly, for $y=0$ we have $0\leq X_{m+1}\leq\theta^{-(n-m)}$ so that the finite-dimensional distributions converge to the deterministic zero path. Thus, in these cases, there is no non-trivial limit in the sense of Theorem~\ref{maintheorem_ge1}. For $a=0$ and $\theta>1$, the same arguments show that in this case we cannot have a non-trivial limit in the sense of Theorem~\ref{maintheorem_ge1} for any starting point.
\end{Remark}
\begin{Remark}
Let $\theta>1$ and $-1<a<0$; this is the gray region     in Figure~\ref{fig:map} and the white region below it. In both parts, every starting point is positive. The persistence restrictions $X_{k+1}>\theta X_k$ force each subsequent variable $X_{k+1}$ to be larger than the previous $X_k$, which is either impossible even in the first step (white region with $-1<a\leq -\bar\theta$) or becomes impossible after finitely many steps (grey region with $-\bar\theta< a<0$). Thus, here the event we want to condition on is empty for large enough $n$ and we cannot expect a  limit in the sense of Theorem~\ref{maintheorem_ge1}.
\end{Remark}
\begin{Remark}
    The white case (W) is trivial: If  $\theta<-1$ and $-1<a\leq0$, then 
    $[-a,1]\subseteq[0,1]$. Hence, $X_i\geq0$ and $\theta X_i\leq0\leq X_{i+1}$ for every $i\in\N$, so every persistence constraint is satisfied automatically. Consequently, $p_n^{a,y}(\theta)=1$, $h(y)=1$ and $\lambda=1$.
    In this case, conditioning on $\Omega_n$ has no effect. 
\end{Remark}

While Theorem~\ref{maintheorem_ge1} determines the limiting conditioned Markov chain, the next main result contains the stationary distribution of the limiting Markov chain.

\begin{Theorem}\label{stationary_distribution}
Let $\Omega$, $\lambda$, and $h$ be as in Theorem~\ref{maintheorem_ge1}. In the light blue region (B), that is, $\theta>1$ and $a>0$, we set 
\begin{align*}
   g(x)\coloneqq  g_{\theta,a}(x)
    \coloneqq
    E\left(\bar\theta,
    \frac{\bar\theta x}{\lambda(1+a)}\right),
     \qquad x\in\Omega=[-a,0).
\end{align*}
In the light green region (G), which is $\theta<-1$ and $a\geq-\bar\theta$, we set
\begin{align*}
    g(x)\coloneqq g_{\theta,a}(x)
    \coloneqq
    E\left(\bar\theta,
    \frac{-\bar\theta x}{\lambda(1+a)}\right),
    \qquad x\in\Omega=(\bar\theta,1].
\end{align*}
In the light yellow region (Y), so $\theta<-1$ and $0<a<-\bar\theta$, we set, on $\Omega=[-a,1]$,
\begin{align*}
    g(x)\coloneqq g_{\theta,a}(x)
    \coloneqq
    \begin{cases}
        E\left(\bar\theta,
        \dfrac{-\bar\theta x}{\lambda(1+a)}\right),
        &-a\leq x\leq-a\theta,\\
        E\left(\bar\theta,
        \dfrac{a}{\lambda(1+a)}\right),
        &-a\theta<x\leq1.
    \end{cases}
\end{align*}
Then $g$ is strictly positive on $\Omega$ and satisfies the left eigenvalue equation
\begin{align}\label{left_eigenvalue_eq}
    \lambda g(y)
    =
    \frac{1}{1+a}
    \int_{\Omega\cap\{y>\theta x\}}
    g(x)\d x,
    \qquad y\in\Omega.
\end{align}
Furthermore, the probability measure
\begin{align*}
   \pi(\d x)\coloneqq \pi_{\theta,a}(\d x)
    \coloneqq
    \frac{\mathds{1}_{\Omega}(x)h(x)g(x)}
    { \int_\Omega
    h(u)g(u)\d u}\d x
\end{align*}
is stationary for the reduced transition kernel $p'$. 
\end{Theorem}
\begin{Remark}
    In the white region, where $\theta<-1$ and $-1<a\leq0$, the stationary distribution is the uniform distribution on $[-a,1]$.
\end{Remark}
Dualities for the adjoint of the killed transition operator connect the expanding and non-expanding
regimes by exchanging right and left eigenfunctions of that operator.
The product formula in Theorem~\ref{stationary_distribution}
therefore also gives stationary distributions for the limiting
chains in the non-expanding regime.
\begin{Corollary}\label{cor_stationary_nonexpanding}
    Let $a>-1$ and $\theta\in[-1,1)$.
   Let $h=h_{\theta,a}$ and $\lambda$ be as in
    \cite[Theorem~1.2]{AW}, with $h_{\theta,a}$ being
    the function for the indicated   parameters restricted to $[-a,1]$.
    Define
    \begin{align*}
    g(x)\coloneqq
    \begin{cases}
        1,&-1\leq\theta<0,\quad -1<a\leq0,\\
        h_{\bar\theta,a}(x),&-1\leq \theta<0,\quad a>0, \\ 
        \mathds{1}_{(0,1]}(x),&\theta=0.
    \end{cases}
\end{align*}
    For $0<\theta<1$, define $g(x)=0$ for $x\leq 0$ and 
    \begin{align*}
        g(x)\coloneqq
            1+\sum_{m=1}^{\infty}
            \frac{\left(\frac{-1}{\theta\lambda(1+a)}\right)^m
            \bar\theta^{m(m-1)/2}}{m!}
            (\theta^m-x)_+^m,\qquad x>0.
    \end{align*}
Then, $g$ satisfies the left eigenvalue equation \eqref{left_eigenvalue_eq} for $\Omega=[-a,1]$ and the probability measure
\begin{align*}
   \pi(\d x)\coloneqq \pi_{\theta,a}(\d x)
    \coloneqq
    \frac{\mathds{1}_{[-a,1]}(x)h(x)g(x)}
    { \int_{[-a,1]}
    h(u)g(u)\d u}\d x
\end{align*}
is stationary for the reduced transition kernel $p'$. 
\end{Corollary}
\begin{Remark}
    In the setting of Corollary~\ref{cor_stationary_nonexpanding},
    the stationary distribution $\pi$ is concentrated on
    \begin{align*}
        [-a,1]\cap(\theta,1],\qquad & -1\leq\theta<0,\\
        [-a,1]\cap(0,1],\qquad &~~~~ ~0\leq\theta<1.
    \end{align*}
    Indeed, the stationary density is proportional to $hg$,
    and $g$ vanishes outside the indicated interval.
\end{Remark}

\paragraph{Related work.}
The closest papers are \cite{AurzadaRaschel,AW}. In \cite{AurzadaRaschel}, the initial starting point has the same uniform distribution as the subsequent
innovations, and the generating functions of the
persistence probabilities are computed in all
parameter regions. The paper \cite{AW} studies persistence with a fixed starting point and the limiting conditioned process in the non-expanding regime $|\theta|\le 1$. The present paper continues this study in the expanding regime.

Persistence for general Markov chains is studied in \cite{AMZ}. This includes moving-average processes and the relation between the persistence rate and the leading eigenvalue of a killed operator. Earlier work on MA(1)-type stationary sequences includes \cite{15}, which uses generating functions. For Gaussian MA(1) processes, perturbation methods are used in \cite{2}, and further estimates for stationary Gaussian processes are given in \cite{12}. The papers \cite{aurzada2025persistenceprobabilitiesma1sequences,AW2} study persistence probabilities and conditioned processes for Laplace innovations. Related questions for autoregressive processes are treated in \cite{1,3,9,DonnartSimon,11,14}. For surveys of persistence probabilities and their role in statistical physics, we refer to \cite{5,8}.

Doob transforms provide a general way to describe conditioned processes; see \cite{doob,RogersWilliams1,RogersWilliams2}. A presentation for Markov chains on countable state spaces is given in \cite{Swart}. Bertoin and Doney \cite{BertoinDoney} study random walks conditioned to stay non-negative. Related conditioning problems arise for ordered random walks \cite{Wachtel7,EichelsbacherKonig}, including their connection with the Airy line ensemble \cite{Wachtel1}. Further limit theorems concern random walks in cones \cite{Wachtel5,Wachtel3} and integrated random walks conditioned to stay positive \cite{Wachtel2}.

Quasi-stationary distributions and processes conditioned on non-absorption are studied in \cite{ChampagnatVillemonais2016,ColletMartinezSanMartin,SenetaVereJones}. Harmonic functions are also used to study stationary measures and Markov chains with asymptotically zero drift; see \cite{Wachtel4,Wachtel6}. In our model, we give explicit formulas for stationary distributions of the limiting chains and quasi-stationary distributions of the restricted killed kernels.

\paragraph{Outline.}
Section~\ref{sec_prelim} recalls the required results from \cite{AW} and establishes the probability and operator dualities. Section~\ref{sec_light_blue} treats $\theta>1$ using generating functions of the persistence probabilities conditioned on a fixed endpoint. Section~\ref{sec_light_green_yellow} treats $\theta<-1$ using spectral methods and derives explicit sum formulas for the eigenfunction. Section~\ref{sec_stationary} identifies the stationary distributions of the limiting chains in the expanding and non-expanding regime. 

\section{Preliminaries}\label{sec_prelim}
\subsection{Auxiliary results}
 Let $\mu_a$ denote the uniform probability measure on $[-a,1]$. For $\theta\in\R$, define the killed transition operator $T_{\theta,a}$ by
\begin{align*}
    T_{\theta,a}f(x)
    \coloneqq\int_{\{y>\theta x\}\cap [-a,1]}f(y)\frac{1}{1+a}\d y,
    \qquad x\in [-a,1],~f\in L^2([-a,1],\mu_a).
\end{align*}
The following auxiliary results are taken from, or are
immediate adaptations of, \cite{AW}.
\begin{Proposition}
\label{spectral_input_contractive}
    Let $\theta\in(-1,1)\setminus\{0\}$ and $a>0$. Denote by
    $\lambda(\theta,a)$ the spectral radius of $T_{\theta,a}$.
    \begin{enumerate}
        \item[(B,G)] If $\theta\in(0,1)$ and $a> 0$ or $\theta\in(-1,0)$ and $a\geq-\theta$, then
        $\lambda(\theta,a)=z_0^{-1}$, where $z_0$ is the smallest
        positive zero of
        \begin{align*}
            z\mapsto E\left(\theta,\frac{-z}{1+a}\right).
        \end{align*}
        \item[(Y)] If $\theta\in(-1,0)$ and $a\in(0,-\theta)$ then
        $\lambda(\theta,a)=z_0^{-1}$, where $z_0$ is the smallest positive zero of
        \begin{align*}
            z\mapsto E\left(\theta,
    \frac{za}{\theta(1+a)}\right)
    -z\frac{\theta+a}{\theta(1+a)}
    E\left(\theta,\frac{za}{1+a}\right).
        \end{align*}
    \end{enumerate}
    In each case, $z_0$ is real, strictly positive and simple. Moreover, for every other zero $z\in\C$, it holds $|z|>z_0$.
\end{Proposition}
\begin{proof}
    The identification of $\lambda$ as the persistence exponent of $p_n^{a,x}(\theta)$ follows from \cite[Theorems~1.1--1.2]{AW}. The identification as the spectral radius of $T_{\theta,a}$ can be seen as follows. The operator describes one step with respect to the persistence
    constraint. In particular,
    \begin{align*}
        T_{\theta,a}^n\mathds{1}(x)=p_n^{a,x}(\theta).
    \end{align*}
    The persistence asymptotics from \cite{AW} therefore
    identify its spectral radius: $T_{\theta,a}$ has a bounded kernel and an underlying finite measure, so it is compact by \cite[Proposition~II.4.7]{ConwayFunctionalAnalysis}. Therefore, we can apply \cite[Theorem~1.2, p.~1413]{AMZ}, to
    $T_{\theta,a}$ on $C([-a,1])$ with initial
    distribution $\mu_a$, which yields that its spectral radius equals
    the exponential persistence rate identified in
    \cite{AW}. 
    The same spectral radius is obtained on
    $L^2([-a,1],\mu_a)$. Indeed, the operator is compact
    on both spaces and maps $L^2(\mu_a)$ continuously
    into $C([-a,1])$. Thus, every eigenfunction for a
    non-zero eigenvalue has a continuous representative,
    and the two operators have the same non-zero spectrum.

    The root formulas are obtained in \cite[Theorems~1.1--1.2]{AW}.
    In cases (B,G), dominance and simplicity
    of the relevant zero follow from \cite[Lemma~2.5]{AW}. In case (Y),
    positivity and uniqueness of the dominant zero follow from
    \cite[Lemma~5.6]{AW}, while simplicity follows from the derivative
    calculation in the proof of \cite[Lemma~5.7]{AW}.
\end{proof}
\begin{Remark}
    The case (B,G) corresponds to the dark blue and dark green region in Figure~\ref{fig:map}; the case (Y) corresponds to the dark yellow region.
\end{Remark}
\begin{Lemma}\label{Prob=Lim_int_trans_kernel}
Assume that there is a continuous function $h:[-a,1]\to[0,\infty)$ which is strictly positive on some set $\Omega\subseteq[-a,1]$ and equal to zero on $[-a,1]\setminus\Omega$, and constants $\lambda=\lambda(\theta,a)>0$, $c(\theta,\lambda,a)>0$ not depending on $x$ or $n$ such that
  $$\lim_{n\to\infty}\frac{p_n^{a,x}(\theta)}{c(\theta,\lambda,a)\lambda^{n+1}}=h(x)$$
  holds for all $x\in[-a,1]$.
Extend $h$ to $\R$ by setting
\begin{align*}
    h(x)\coloneqq\frac{1}{\lambda(1+a)}
    \int_{\{y>\theta x\}\cap\Omega}h(y)\d y,
    \qquad x\in\R\setminus[-a,1].
\end{align*}
Then, this extension is continuous and satisfies~\eqref{eq_eigenfunction}. Moreover, for any $y\in \Omega$, the limit in \eqref{def_lim_process}
exists, and $M'$ is a Markov chain with reduced transition kernel 
\begin{align*} 
    p'(x_1,\d x_2)\coloneqq \mathds{1}_S(x_2,x_1)\frac{h(x_2)}{h(x_{1})}\frac{1}{\lambda}\frac{1}{a+1}\d x_2,~~~~~~~~x_1,x_2\in \Omega.
\end{align*}
\end{Lemma}
\begin{proof}
For $x\notin[-a,1]$, the eigenvalue equation holds by the definition of the extension. On $[-a,1]$, the proof of \cite[Lemma~2.2]{AW} gives the eigenvalue equation from the assumed asymptotic relation. Its right-hand side is continuous in $x$, so the extension is continuous on $\R$. 

The argument for the limit in \eqref{def_lim_process} is the same as in \cite[Lemma~2.3]{AW}. The required domination follows from $p_n^{a,x}(\theta)\leq p_n^{a,-a}(\theta)$ for $\theta\geq0$ and $p_n^{a,x}(\theta)\leq p_n^{a,1}(\theta)$ for $\theta\leq0$. The proof remains valid when $h$ vanishes outside $\Omega$: the initial point belongs to $\Omega$, so its limiting denominator is positive, while a path whose last coordinate lies outside $\Omega$ receives limiting weight zero. Finally, the eigenvalue equation \eqref{eq_eigenfunction} shows that $p'$ has total mass one.
\end{proof}

\subsection{Dualities}
The key observation is that time reversal turns a fixed starting point into a fixed endpoint. We first state this relation for persistence probabilities and then for killed transition operators. In this way, we reduce the expanding parameter $\theta$ to the contractive parameter $\bar\theta=\theta^{-1}$ and denote the endpoint-fixed quantities by $r_0^{a,y}(\theta)\coloneqq1$ and
\begin{align*}
    r_n^{a,y}(\theta)\coloneqq\P(X_2\geq\theta X_1,X_3\geq\theta X_2,\dots,y\geq\theta X_n),\quad n\ge 1.
\end{align*}
 We also introduce the generating functions
 \begin{align*}
    \hat P_{a,y,\theta}(z)\coloneqq\sum_{n=0}^{\infty}p_n^{a,y}(\theta)z^n\quad\text{and}\quad\hat R_{a,y,\theta}(z)
    \coloneqq\sum_{n=0}^{\infty}r_n^{a,y}(\theta)z^n.
\end{align*}
\begin{Lemma}
    Let $\theta>1$, $a>0$ and $\bar\theta=\theta^{-1}$. Then, for all $n\in\N_0$,
       $p_n^{a,y}(\theta)
        =r_n^{1/a,-y/a}(\bar\theta)$
    and
    \begin{align}
       \hat P_{a,y,\theta}(z)
        &=\hat R_{1/a,-y/a,\bar\theta}(z).\label{fixed_start_endpoint_duality_generating_functions_theta_ge1}
    \end{align}
\end{Lemma}
\begin{proof}
The proof follows by dividing the persistence restrictions by $-a$, using that $\tilde X_i\coloneqq\frac{X_i}{-a}\sim U\left[-1/a,1\right]$ and $(\tilde X_2,\dots,\tilde X_{n+1})\overset{d}{=}(\tilde X_n,\dots,\tilde X_1)$.
        Summation of the probability identity in $n$ proves the claim for the generating function.
\end{proof}
Next, we formulate these dualities in terms of operators. By definition and the Markov property, we have, for $y\in[-a,1]$,
\begin{align*}
    p_n^{a,y}(\theta)
    =\frac{1}{1+a} \int_{[-a,1]\cap\{x>\theta y\}}p_{n-1}^{a,x}(\theta)\d x
    =T_{\theta,a}^n\mathds{1}(y),
    \qquad n\ge 1.
\end{align*}
The adjoint operator with respect to $L^2([-a,1],\mu_a)$ is given by
\begin{align*}
    T_{\theta,a}^*f(y)
    =\int_{\{y>\theta x\}\cap [-a,1]}f(x)\frac{1}{1+a}\d x,
    \qquad y\in [-a,1].
\end{align*}
The two operators use the same constraint in opposite directions: in $T_{\theta,a}$, the variable $x$ is fixed and $y$ is integrated out, whereas in $T_{\theta,a}^*$, the endpoint $y$ is fixed and the possible preceding values $x$ are integrated out. In particular, we have
\begin{align*}
    r_n^{a,y}(\theta)
    =\frac{1}{1+a}\int_{[-a,1]\cap\{y>\theta x\}} r_{n-1}^{a,x}(\theta)\d x
    =\left(T_{\theta,a}^*\right)^n\mathds{1}(y),\qquad n\geq1.
\end{align*}

This gives the following operator identities.
\begin{Lemma}
    Let $\theta\neq 0$ and recall $\bar\theta=\theta^{-1}$.
    \begin{enumerate}[label={\rm (\roman{*})}]
        \item If $\theta<0$, then
        \begin{align}\label{adjoint_duality_negative_theta}
            T_{\theta,a}=T_{\bar\theta,a}^*.
        \end{align}
        \item If $\theta>0$ and $a>0$, define the map $U_a:L^2([-a,1],\mu_a)\to L^2([-1/a,1],\mu_{1/a})$ by
        \begin{align*}
            U_af(u)\coloneqq f(-au),
            \qquad u\in [-1/a,1].
        \end{align*}
        Then, $U_a$ is unitary and
        \begin{align}\label{adjoint_duality_positive_theta}
            U_aT_{\theta,a}U_a^{-1}=T_{\bar\theta,1/a}^*.
        \end{align}
    \end{enumerate}
\end{Lemma}
\begin{proof}
   The claim for $\theta<0$ follows since for $x,y\in[-a,1]$, the inequality $y>\theta x$ is equivalent to $x>\bar\theta y$.     Now let $\theta>0$ and $a>0$. The map $U_a$ is unitary, since it is bijective and
\begin{align*}
    \int_{-1/a}^1 |f(-au)|^2 \mu_{1/a}(\d u)
    =
    \int_{-a}^1 |f(x)|^2 \mu_a(\d x).
\end{align*}
For a function $f$ on $[-1/a,1]$, we obtain the claim:
    \begin{align*}
        U_aT_{\theta,a}U_a^{-1}f(u)
        &=\int_{[-a,1]}\mathds{1}_{\{y>-a\theta u\}}f\left(\frac{-y}{a}\right)\frac{1}{1+a}\d y
        =\int_{[-1/a,1]}\mathds{1}_{\{-av>-a\theta u\}}f(v)\frac{a}{1+a}\d v
        \\&=\int_{[-1/a,1]}\mathds{1}_{\{u>\bar\theta v\}}f(v)\mu_{1/a}(\d v)
        =T_{\bar\theta,1/a}^*f(u). \qedhere
    \end{align*}
\end{proof}

\section{The light blue region}\label{sec_light_blue}
We first compute the endpoint-fixed generating function for a parameter $\theta\in (0,1)$. The dualities then give the startpoint-fixed generating function for a parameter $\theta\in(1,\infty)$. This is the adjoint version of the generating-function calculation in \cite{AW}.

The proof proceeds as follows. Repeated multiplication by $\theta\in(0,1)$ produces nested intervals whose endpoints converge to zero. The positive-part notation $x_+=\max(x,0)$ for $x\in\R$ combines the resulting expressions on the respective intervals into a single formula.

\begin{Lemma}\label{endpoint_gf_positive_contractive_parameter}
    Let $\theta\in(0,1)$ and $b>0$.
    For $y\in[-b,0]$, define
    \begin{align*}
        R_-(y,z)
        \coloneqq1+\sum_{m=1}^{\infty}
        \frac{\left(\frac{z}{\theta(1+b)}\right)^m\theta^{-m(m-1)/2}}{m!}
        \left(y+b\theta^m\right)_+^m,
    \end{align*}
    and for $y\in[0,1]$, define
    \begin{align*}
        R_+(y,z)
        \coloneqq1+\sum_{m=1}^{\infty}
        \frac{\left(\frac{-z}{\theta(1+b)}\right)^m\theta^{-m(m-1)/2}}{m!}
        \left(\theta^m-y\right)_+^m.
    \end{align*}
    Then, within the radius of convergence of the power series on the left-hand side, we have
    \begin{align}\label{formula_endpoint_gf_positive_contractive_parameter}
        \hat R_{b,y,\theta}(z)
        =\begin{cases}
            R_-(y,z),&-b\leq y\leq0,
            \\\displaystyle
            \frac{E\left(\theta,\frac{zb}{1+b}\right)}
            {E\left(\theta,\frac{-z}{1+b}\right)}R_+(y,z),&0< y\leq1.
        \end{cases}
    \end{align}
\end{Lemma}
\begin{proof}
   We first work with $|z|$ sufficiently small. Since $0\leq r_n^{b,y}(\theta)\leq1$, the defining series $\hat R_{b,y,\theta}(z)$ converges uniformly in $y\in[-b,1]$ whenever $|z|<1$. Since $r_n^{b,\cdot}(\theta)=T_{\theta,b}^*r_{n-1}^{b,\cdot}(\theta)$, for $n\geq1$, we obtain
    \begin{align}\label{integral_equation_endpoint_gf_positive}
        \hat R_{b,y,\theta}(z)
        &=1+\sum_{n=1}^{\infty}r_n^{b,y}(\theta)z^n\notag
        =1+\sum_{n=1}^{\infty}
        \int_{-b}^{\min(y\bar\theta,1)}
        r_{n-1}^{b,x}(\theta)z^n\frac{1}{1+b}\d x\notag\\
        &=1+\frac{z}{1+b}
        \int_{-b}^{\min(y\bar\theta,1)}
        \sum_{n=0}^{\infty}r_n^{b,x}(\theta)z^n\d x
        =1+\frac{z}{1+b}
        \int_{-b}^{\min(y\bar\theta,1)}
        \hat R_{b,x,\theta}(z)\d x.
    \end{align}
    Here, the integral is understood to be zero if $y\bar\theta\leq-b$.
    In the third equality, we shifted the summation index and
    interchanged the sum and integral, which is justified by uniform
    convergence for $|z|<1$.
      Consequently, the end-point generating function satisfies the differential equation
    \begin{align}\label{differential_equation_endpoint_gf_positive}
        \frac{\partial}{\partial y}\hat R_{b,y,\theta}(z)
        =\frac{z}{\theta(1+b)}\hat R_{b,y\bar\theta,\theta}(z),
        \qquad -b\theta<y<\theta.
    \end{align}
    If $-b\leq y\leq-b\theta$, then $y\bar\theta\leq-b$, so the integral is empty and $\hat R_{b,y,\theta}(z)=1$. If $\theta\leq y\leq1$, then $y\bar\theta\geq1$, so the upper integration boundary is fixed at one and the function is constant in $y$. These two boundary values determine the solution step by step towards zero.

    We now check that the function in \eqref{formula_endpoint_gf_positive_contractive_parameter} satisfies the differential equation, the boundary conditions and continuity at zero. We then verify the integral equation and uniqueness.

    For fixed $z$, the function $y\mapsto R_-(y,z)$ is continuously
    differentiable on $(-b\theta,0)$, including the points
    $y=-b\theta^m$, $m\geq2$. Indeed, the series is locally
    finite, and $t\mapsto t_+^m$ is continuously differentiable
    for $m\geq2$. Differentiating the defining series for $R_-$ and shifting the index from $m$ to $m-1$ gives, for $ -b\theta<y<0$,
    \begin{align*}
        \frac{\partial}{\partial y}R_-(y,z)
         &=\sum_{m=1}^{\infty}\frac{\left(\frac{z}{\theta(1+b)}\right)^m\theta^{-m(m-1)/2}}{(m-1)!}(y+b\theta^m)_+^{m-1}\\
         &=\frac{z}{\theta(1+b)}\sum_{m=0}^{\infty}\frac{\left(\frac{z}{\theta(1+b)}\right)^m\theta^{-m(m+1)/2}}{m!}(y+b\theta^{m+1})_+^m\\
        &=\frac{z}{\theta(1+b)}\sum_{m=0}^{\infty}\frac{\left(\frac{z}{\theta(1+b)}\right)^m\theta^{-m(m-1)/2}}{m!}\left(\frac{y}{\theta}+b\theta^m\right)_+^m
        =\frac{z}{\theta(1+b)}R_-\left(\frac{y}{\theta},z\right), 
    \end{align*}
    where we used $(y+b\theta^{m+1})_+^m=\theta^m(y/\theta+b\theta^m)_+^m$.
    Similarly, $R_+(\cdot,z)$ is continuously differentiable
    on $(0,\theta)$, including the points $y=\theta^m$,
    $m\geq2$. The same calculation for $R_+$, together with the minus sign obtained by differentiating $(\theta^m-y)_+^m$, gives
    \begin{align*}
        \frac{\partial}{\partial y}R_+(y,z)
        &=\frac{z}{\theta(1+b)}R_+\left(\frac{y}{\theta},z\right), \qquad 0<y<\theta.
    \end{align*}
    Thus, both functions satisfy \eqref{differential_equation_endpoint_gf_positive} on their respective domains. Moreover, $R_-(y,z)=1$ for $-b\leq y\leq-b\theta$ and $R_+(y,z)=1$ for $\theta\leq y\leq1$. Therefore, if we multiply $R_+(y,z)$ with the constant ${E\left(\theta,\frac{zb}{1+b}\right)}/
    {E\left(\theta,\frac{-z}{1+b}\right)}$, we get the boundary pieces explained above. 
Dominated convergence 
allows us to take the corresponding one-sided limits at zero
term by term, which gives
\begin{align*}
    \lim_{y\searrow0}
    \frac{E\left(\theta,\frac{zb}{1+b}\right)}
         {E\left(\theta,\frac{-z}{1+b}\right)} R_+(y,z)
    &=\frac{E\left(\theta,\frac{zb}{1+b}\right)}{E\left(\theta,\frac{-z}{1+b}\right)}
    \sum_{m=0}^{\infty}
    \frac{\theta^{m(m-1)/2}}{m!}
    \left(\frac{-z}{1+b}\right)^m
    \\&=E\left(\theta,\frac{zb}{1+b}\right)
    \\&=R_-(0,z).
\end{align*}
Thus, the two branches in
\eqref{formula_endpoint_gf_positive_contractive_parameter}
define a continuous function at zero.

It remains to verify
\eqref{integral_equation_endpoint_gf_positive}.
Let $R(y,z)$ denote the proposed function on the right-hand side
of \eqref{formula_endpoint_gf_positive_contractive_parameter}.
For $y\in(-b\theta,0)\cup(0,\theta)$, the differential equation \eqref{differential_equation_endpoint_gf_positive} gives
\begin{align*}
    &\frac{\partial}{\partial y}
    \left[
        R(y,z)-1-\frac{z}{1+b}
        \int_{-b}^{y/\theta}R(x,z)\d x
    \right]
    =\frac{\partial}{\partial y}R(y,z)-\frac{z}{\theta(1+b)}R(y/\theta,z)
    =0.
\end{align*}
Hence, the difference in brackets is constant on each of
$(-b\theta,0)$ and $(0,\theta)$. At the left boundary, we have
\begin{align*}
    R(-b\theta,z)-1
    -\frac{z}{1+b}\int_{-b}^{-b}R(x,z)\d x
    =1-1=0.
\end{align*}
By continuity, the difference is therefore zero on
$[-b\theta,0)$. Its continuity at zero then shows that it is
also zero on $[0,\theta]$.
For $y\in[-b,-b\theta]$, the integral is empty and $R(y,z)=1$.
For $y\in[\theta,1]$, both $R(y,z)$ and the upper integration
boundary are constant, so the identity follows from its value
at $y=\theta$. This proves the integral equation on $[-b,1]$.

It remains to prove uniqueness. Since $T_{\theta,b}^*$ has operator norm at most one
on $C([-b,1])$, the operator $I-zT_{\theta,b}^*$
is invertible for $|z|<1$ by
\cite[Lemma~VII.2.1, p.~192]{ConwayFunctionalAnalysis}, where $I$ is the identity operator.
Hence, the integral equation has the unique continuous solution
\begin{align*}
    f=(I-zT_{\theta,b}^*)^{-1}\mathds{1}
\end{align*} 
and the proposed function agrees with the
generating function for $|z|$ sufficiently small.
The identity theorem extends the formula throughout the disk
of convergence of the defining power series, with removable
singularities on the right-hand side filled in by continuity.
\end{proof}
We now return to $\theta>1$. By \eqref{fixed_start_endpoint_duality_generating_functions_theta_ge1},
we apply Lemma~\ref{endpoint_gf_positive_contractive_parameter}
with parameter $\bar\theta$, with $b=1/a$, and with endpoint
$u=-y/a$.
The change of variables switches the two halves of the interval: $y<0$ corresponds to $u>0$ and therefore to the second line of \eqref{formula_endpoint_gf_positive_contractive_parameter}, whereas $y>0$ corresponds to $u<0$ and to its first line. This explains why only the expression for negative starting points in the following corollary contains a denominator.
\begin{Corollary}
    Let $\theta>1$, $a>0$ and $\bar\theta=\theta^{-1}$. Then, for $y\in[-a,1]$,
    \begin{align}\label{explicit_generating_function_theta_greater_one}
        \hat P_{a,y,\theta}(z)
        =\begin{cases}
        \displaystyle
        \frac{E\left(\bar\theta,\frac{z}{1+a}\right)}
        {E\left(\bar\theta,\frac{-az}{1+a}\right)}
        \left(
        1+\sum_{m=1}^{\infty}
        \frac{\left(\frac{-\theta z}{1+a}\right)^m
        \theta^{m(m-1)/2}}{m!}
        \left(a\bar\theta^m+y\right)_+^m
        \right),&-a\leq y<0,
        \\\displaystyle
        1+\sum_{m=1}^{\infty}
        \frac{\left(\frac{\theta z}{1+a}\right)^m
        \theta^{m(m-1)/2}}{m!}
        \left(\bar\theta^m-y\right)_+^m,&0\leq y\leq1.
        \end{cases}
    \end{align}
\end{Corollary}
\begin{proof}
    For $|z|$ sufficiently small, the formula follows from \eqref{fixed_start_endpoint_duality_generating_functions_theta_ge1} and Lemma~\ref{endpoint_gf_positive_contractive_parameter} with $\frac{1}{1+1/a}=\frac{a}{1+a}$ and $u/\bar\theta=-\frac{\theta y}{a}$.
    After inserting $u=-y/a$, the factor $a^{-m}$ from $(\bar\theta^m-u)_+^m$ cancels the factor $a^m$ produced by $1/(1+1/a)^m$, which gives the first line of \eqref{explicit_generating_function_theta_greater_one}. The other branch is obtained analogously. The identity holds within the disk of convergence of the
generating function on the left-hand side.
\end{proof}
We can now identify the dominant pole and the corresponding eigenfunction. By Remark~\ref{Case_B_y_pos}, persistence is possible for only
finitely many steps when $y>0$, while its probability decays
faster than exponentially when $y=0$. Thus, the generating
function is a polynomial in the first case and entire in the
second. We therefore consider $y<0$. In this case, the only denominator depending on $z$ in \eqref{explicit_generating_function_theta_greater_one} is $E(\bar\theta,-az/(1+a))$. Let $z_0$ be its smallest positive zero, so
\begin{align}\label{root_theta_greater_one}
    E\left(\bar\theta,\frac{-az_0}{1+a}\right)=0
\end{align}
and set $\lambda\coloneqq z_0^{-1}$. By \eqref{adjoint_duality_positive_theta} and Proposition~\ref{spectral_input_contractive}, the spectral radius of $T_{\theta,a}$ as well as $T_{\bar\theta,1/a}$ is given by $\lambda$. Moreover, Proposition~\ref{spectral_input_contractive} shows that
$z_0$ is the unique zero of the denominator in the first line of
\eqref{explicit_generating_function_theta_greater_one} of smallest
modulus. Hence, the right-hand side of
\eqref{explicit_generating_function_theta_greater_one} is analytic on
$|z|<z_0$. Since it agrees with the defining power series near zero,
the identity holds throughout this disc. For $y<0$, let $B_y(z)$ denote the finite-sum factor in parentheses in the first line of \eqref{explicit_generating_function_theta_greater_one}. 
Define $h:\R\to\R$ by
\begin{align}\label{eigenfunction_theta_greater_one}
    h(y)
    \coloneqq\begin{cases}
    \displaystyle
    1+\sum_{m=1}^{\infty}
    \frac{\left(\frac{-\theta}{\lambda(1+a)}\right)^m
    \bar\theta^{-m(m-1)/2}}{m!}
    \left(a\bar\theta^m+y\right)_+^m,&y<0,
    \\0,&y\geq0.
    \end{cases}
\end{align}
Thus, we have $h(y)=B_y(z_0)$ for $ -a\leq y<0$.
For each fixed $y<0$, the sum in \eqref{eigenfunction_theta_greater_one} is finite. We next prove continuity and strict positivity of $h$ on $(-\infty,0)$.
\begin{Lemma}\label{positive_eigenfunction_theta_greater_one}
    Let $\theta>1$ and $a>0$. The function $h$ defined in \eqref{eigenfunction_theta_greater_one} is continuous and non-negative on $\R$, and strictly positive on $\Omega=[-a,0)$.
    Moreover, $h(y)=1$ for $y\leq-a\bar\theta$ and $h(y)=0$ for $y\geq0$.
\end{Lemma}
\begin{proof}
    We first prove continuity. If $y\leq-a\bar\theta$, then $a\bar\theta^m+y\leq0$ for every $m\geq1$. Hence, all summands in \eqref{eigenfunction_theta_greater_one} vanish and $h(y)=1$. By definition, $h(y)=0$ for $y\geq0$.
    On $(-\infty,0)$, the defining series is locally finite.
    Moreover, for $y\leq0$ and $m\geq1$,
    \begin{align*}
        \left|
        \frac{\left(\frac{-\theta}{\lambda(1+a)}\right)^m
        \bar\theta^{-m(m-1)/2}}{m!}
        \left(a\bar\theta^m+y\right)_+^m
        \right|
        &\leq
        \frac{\left(\frac{\theta}{\lambda(1+a)}\right)^m
        \bar\theta^{-m(m-1)/2}}{m!}
        \left(a\bar\theta^m\right)^m
        \\&=
        \frac{\bar\theta^{m(m-1)/2}}{m!}
        \left(\frac{a}{\lambda(1+a)}\right)^m.
    \end{align*}
    The bound is independent of $y$, and
    \begin{align*}
        \sum_{m=0}^{\infty}
        \frac{\bar\theta^{m(m-1)/2}}{m!}
        \left(\frac{a}{\lambda(1+a)}\right)^m
        =
        E\left(\bar\theta,\frac{a}{\lambda(1+a)}\right)
        <\infty.
    \end{align*}
    Dominated convergence therefore gives
    \begin{align*}
        \lim_{y\nearrow0}h(y)=\sum_{m=0}^{\infty}\frac{\bar\theta^{m(m-1)/2}}{m!}\left(\frac{-a}{\lambda(1+a)}\right)^m=E\left(\bar\theta,\frac{-a}{\lambda(1+a)}\right)=0,
    \end{align*}
    where the last equality follows from \eqref{root_theta_greater_one}. Thus, $h$ is continuous on $\R$.

    For the positivity, we first show that $h$ is non-negative. Fix $y\in[-a,0)$. For every real $z\in[0,z_0)$,
        $E\left(\bar\theta,\frac{z}{1+a}\right)>0$
    because all coefficients of this power series are positive. Moreover,
        $E\left(\bar\theta,\frac{-az}{1+a}\right)>0,$
    since this function is equal to one at zero and $z_0$ is its smallest positive zero. Further, $\hat P_{a,y,\theta}(z)>0$ since the coefficients of the series are non-negative and $p_0^{a,y}(\theta)=1$. Hence, the first line of \eqref{explicit_generating_function_theta_greater_one} gives
    \begin{align*}
        B_y(z)
        =\hat P_{a,y,\theta}(z)\frac{E\left(\bar\theta,\frac{-az}{1+a}\right)}{E\left(\bar\theta,\frac{z}{1+a}\right)}
        >0.
    \end{align*}
    Since $B_y$ is continuous,
    \begin{align*}
        h(y)=B_y(z_0)=\lim_{z\nearrow z_0}B_y(z)\geq0.
    \end{align*}
    Together with the constant parts considered above, this shows that $h$ is non-negative on $\R$.

    For $-a\bar\theta<y<0$, termwise differentiation gives
    \begin{align*}
        h'(y)
        &=
        \sum_{m=1}^{\infty}
        \frac{\left(\frac{-\theta}{\lambda(1+a)}\right)^m
        \theta^{m(m-1)/2}}{(m-1)!}
        \left(a\bar\theta^m+y\right)_+^{m-1}\\
        &=
        \frac{-\theta}{\lambda(1+a)}
        \sum_{k=0}^{\infty}
        \frac{\left(\frac{-\theta}{\lambda(1+a)}\right)^k
        \theta^{k(k+1)/2}}{k!}
        \left(a\bar\theta^{k+1}+y\right)_+^k\\
        &=
        \frac{-\theta}{\lambda(1+a)}
        \sum_{k=0}^{\infty}
        \frac{\left(\frac{-\theta}{\lambda(1+a)}\right)^k
        \theta^{k(k-1)/2}}{k!}
        \left(a\bar\theta^k+\theta y\right)_+^k
        =
        \frac{-\theta}{\lambda(1+a)}h(\theta y),
    \end{align*}
    where we used $\left(a\bar\theta^m+\theta y\right)_+^m=\theta^m\left(a\bar\theta^{m+1}+y\right)_+^m$.
    It remains to prove strict positivity. Suppose that $h$
    vanishes at some point of $[-a,0)$, and let $y_0$ be its
    leftmost zero in $[-a,0]$. Since $h=1$ on
    $[-a,-a\bar\theta]$, we have $-a\bar\theta<y_0<0$.
    Moreover, $-a<\theta y_0<y_0$, so $h(\theta y_0)>0$
    by the choice of $y_0$. The differential equation gives
    \begin{align*}
        h'(y_0)
        =-\frac{\theta}{\lambda(1+a)}h(\theta y_0)<0.
    \end{align*}
    This contradicts $h\geq0$, since $h$ has a local minimum
    at $y_0$. Hence, $h>0$ on $[-a,0)$.
\end{proof}

\begin{Lemma}\label{asymptotic_theta_greater_one}
    Let $\theta>1$, $a>0$ and $y\in[-a,0)$. Then, setting $c(\theta,\lambda,a)\coloneqq\frac{1+a}{a}\frac{E\left(\bar\theta,\frac{1}{\lambda(1+a)}\right)}{E\left(\bar\theta,\frac{-a\bar\theta}{\lambda(1+a)}\right)}$,
    \begin{align*}
        p_n^{a,y}(\theta)\sim c(\theta,\lambda,a)h(y)\lambda^{n+1}.
    \end{align*}
    \end{Lemma}
\begin{proof}
    For $y<0$, write the first line of \eqref{explicit_generating_function_theta_greater_one} as $\hat P_{a,y,\theta}(z)=\frac{A(z)B_y(z)}{D(z)}$, where $B_y(z_0)=h(y)$ as given in \eqref{eigenfunction_theta_greater_one} and
    \begin{align*}
        A(z)\coloneqq E\left(\bar\theta,\frac{z}{1+a}\right),\qquad D(z)\coloneqq E\left(\bar\theta,\frac{-az}{1+a}\right).
    \end{align*}
    By Proposition~\ref{spectral_input_contractive} (B,G), the zero $z_0=\lambda^{-1}$ of $D$
    with smallest modulus is unique and simple. Hence, there exists
    $z_1>z_0$ such that $D$ has no other zero in
    $\{z\in\C:|z|<z_1\}$. Moreover, $A(z_0)>0$ since $z_0>0$ and the coefficients of the defining series of $A$ are non-negative and the first one is equal to one. By
    Lemma~\ref{positive_eigenfunction_theta_greater_one}, $B_y(z_0)=h(y)>0$.
    Thus, there is no cancellation at $z_0$, and $z_0$ is the unique
    singularity of $\hat P_{a,y,\theta}$ of smallest modulus. Define, for $|z|<z_1$,
    \begin{align*}
    R_y(z)\coloneqq
    \begin{cases}
    \frac{A(z)B_y(z)}{D(z)}-\frac{A(z_0)h(y)}{D'(z_0)(z-z_0)},&z\neq z_0,
    \\\frac{(A B_y)'(z_0)}{D'(z_0)}-\frac{A(z_0)h(y)D''(z_0)}{2D'(z_0)^2},&z=z_0.
    \end{cases}
    \end{align*}
    For $z\neq z_0$ sufficiently close to $z_0$,  we use Taylor expansion of the $A(z)B_y(z)$ and $D(z)$, respectively, at $z_0$. With $D(z_0)=0$, $D'(z_0)\neq0$ and $B_y(z_0)=h(y)$, this is
    \begin{align*}
        \frac{A(z)B_y(z)}{D(z)}
        &=\frac{ A(z_0)h(y)+(AB_y)'(z_0)(z-z_0)+O((z-z_0)^2)}{(z-z_0)\left(D'(z_0)+\frac12D''(z_0)(z-z_0)+O((z-z_0)^2)\right)}\\
          &=\frac{A(z_0)h(y)+(AB_y)'(z_0)(z-z_0)+O((z-z_0)^2)}{D'(z_0)(z-z_0)}\\&\qquad\cdot\left(1+\frac{D''(z_0)}{2D'(z_0)}(z-z_0)+O((z-z_0)^2)\right)^{-1}\\
        &=\left(\frac{A(z_0)h(y)}{D'(z_0)(z-z_0)}+\frac{(AB_y)'(z_0)}{D'(z_0)}+O(z-z_0)\right)\\
    &\qquad\cdot\left(1-\frac{D''(z_0)}{2D'(z_0)}(z-z_0)+O((z-z_0)^2)\right)\\
        &=\frac{A(z_0)h(y)}{D'(z_0)(z-z_0)}
        +\frac{(AB_y)'(z_0)}{D'(z_0)}
        -\frac{A(z_0)h(y)D''(z_0)}{2D'(z_0)^2}
        +O(z-z_0).
    \end{align*}
    Here, we used $(1+w)^{-1}=1-w+O(w^2)$.
    Thus, the prescribed value $R_y(z_0)$ is precisely
    the limit of $R_y(z)$ as $z\to z_0$.
    The singularity of $R_y$ at $z_0$ is removable, and hence $R_y$ is analytic on $\{z\in\mathbb{C}:|z|<z_1\}$. Therefore,
    \begin{align*}
    \hat P_{a,y,\theta}(z)
    =\frac{A(z_0)h(y)}{D'(z_0)(z-z_0)}+R_y(z)
    =\frac{-A(z_0)h(y)}{z_0D'(z_0)}\frac{1}{1-z/z_0}+R_y(z).
    \end{align*}
   To apply Lemma~2.6 from \cite{AW},
   set
    \begin{align*}
        f(z)\coloneqq \frac{-A(z_0)h(y)}{D'(z_0)}+(z_0-z)R_y(z).
    \end{align*}
    Then, $f$ is analytic on $\{z\in\mathbb{C}:|z|<z_1\}$.  Furthermore, $D'(z_0)=\frac{-a}{1+a}E\left(\bar\theta,\frac{-a\bar\theta z_0}{1+a}\right)=\frac{-a}{1+a}D(\bar\theta z_0)<0$ since $D(0)>0$, $0<\bar\theta z_0<z_0$ and $z_0$ is the first positive zero of
$D$, so  
    \begin{align*}
        \frac{-A(z_0)}{D'(z_0)}=\frac{1+a}{a}\frac{E\left(\bar\theta,\frac{1}{\lambda(1+a)}\right)}{E\left(\bar\theta,\frac{-a\bar\theta}{\lambda(1+a)}\right)}=c(\theta,\lambda,a)>0,
    \end{align*}
    which proves $f(z_0)\ne 0$. Now, Lemma~2.6 from \cite{AW}
    yields
    \begin{align*}
        p_n^{a,y}(\theta)\sim\frac{-A(z_0)}{z_0D'(z_0)}h(y)\lambda^n=\frac{-A(z_0)}{D'(z_0)}h(y)\lambda^{n+1}.\tag*{\qedhere}
    \end{align*}
\end{proof}
\begin{proof}[Proof of Theorem~\ref{maintheorem_ge1} (B)]
    By Lemma~\ref{positive_eigenfunction_theta_greater_one}, the function $h$ is strictly positive on $[-a,0)$. Moreover, Lemma~\ref{asymptotic_theta_greater_one} gives
    \begin{align*}
        p_n^{a,y}(\theta)\sim c(\theta,\lambda,a)h(y)\lambda^{n+1},\qquad y\in[-a,0).
    \end{align*}
   Since $h=0$ on $[0,1]$, Remark~\ref{Case_B_y_pos} gives
   \begin{align*}
       \lim_{n\to\infty}\frac{p_n^{a,y}(\theta)}{c(\theta,\lambda,a)\lambda^{n+1}}=h(y)\qquad y\in[0,1].
   \end{align*}
   Thus, all assumptions of Lemma~\ref{Prob=Lim_int_trans_kernel} are satisfied. The lemma shows that $h$ satisfies the eigenvalue equation, that the limiting finite-dimensional distributions exist and that their reduced transition kernel is given as in Theorem~\ref{maintheorem_ge1} (B). The extension in Lemma~\ref{Prob=Lim_int_trans_kernel} agrees with the stated function outside $[-a,1]$: for $x<-a$, its defining integral is the same as at $x=-a$, so its value is one; for $x>1$, the integration set is empty, so its value is zero.
\end{proof}
\section{The light green and light yellow regions}\label{sec_light_green_yellow}

Let $\theta<-1$ and recall $\bar\theta= \theta^{-1}\in(-1,0)$. By \eqref{adjoint_duality_negative_theta}, $T_{\theta,a}$ is the adjoint of the already treated operator with parameter $\bar\theta$, so no further reflection or rescaling is needed. 

The proof has three steps. First, we use adjoint duality and the results from Proposition~\ref{spectral_input_contractive} to identify the spectral radius $\lambda$. Second, by applying results from operator theory, we prove that $\lambda$ is the only spectral value of modulus $\lambda$ and it has a (up to scaling) unique positive eigenfunction. The eigenvalue equation then gives the polynomial pieces of $h$ and the explicit sum formulas in Lemma~\ref{method_of_steps_theta_smaller_minus_one_operator}. Third, a Riesz projection gives the persistence asymptotics.

First, we introduce some notation. For a bounded linear operator $A$ on a complex Hilbert space $\mathcal H$ with identity operator $I$ on $\mathcal H$, we denote its spectrum by 
\begin{align*}
    \sigma(A)\coloneqq \{z\in\C:zI-A\text{ is not invertible}\}
\end{align*}
and its spectral radius by 
\begin{align*}
    r(A)\coloneqq \sup\{|z|:z\in\sigma(A)\}.
\end{align*}
For an eigenvalue $\nu$, its generalized eigenspace is
\begin{align*}
    \bigcup_{k\geq1}\ker(A-\nu I)^k.
\end{align*}
The eigenvalue is called algebraically simple if this space is one-dimensional. An operator on $L^2$ is called positivity improving if it maps every non-zero non-negative function to a function which is strictly positive
almost everywhere.
Moreover, we define $\operatorname{ran}A\coloneqq\{Af:f\in \mathcal H\}$. For $B\subseteq\C$, we write
    $\overline{B}\coloneqq\{\overline z:z\in B\}.$
The adjoint $A^*$ is characterized by
    $\langle Af,g\rangle=\langle f,A^*g\rangle$, $f,g\in \mathcal H$.

Now we identify the spectral radius. Assume that $a>0$. We choose the characteristic equation according to the position of $(\bar\theta,a)$ in the parameter diagram. If $a\geq-\bar\theta$, the pair lies in the case (B,G) of Proposition~\ref{spectral_input_contractive}, and we let $z_0$ be the smallest positive zero of
\begin{align*}
   z\mapsto E\left(\bar\theta,\frac{-z}{1+a}\right).
\end{align*}
If $0<a<-\bar\theta$, the pair $(\bar\theta,a)$ lies in the case (Y) of Proposition~\ref{spectral_input_contractive}, and we let $z_0$ be the smallest positive zero of 
\begin{align*}
   z\mapsto E\left(\bar\theta,
    \frac{za}{\bar\theta(1+a)}\right)
    -z\frac{\bar\theta+a}{\bar\theta(1+a)}
    E\left(\bar\theta,\frac{za}{1+a}\right).
\end{align*}
In both parameter ranges, Proposition~\ref{spectral_input_contractive} shows that $z_0$ is real, strictly positive and simple and is the unique zero of smallest
modulus. Set $\lambda\coloneqq z_0^{-1}$.
\begin{Lemma}
    The spectral radius of
$T_{\theta,a}$ as well as $T_{\bar\theta,a}$ is given by $\lambda$.
\end{Lemma}
\begin{proof}
    Since an operator is invertible if and only if its adjoint is invertible, we have, by \eqref{adjoint_duality_negative_theta},
$\sigma(T_{\theta,a})=\overline{\sigma(T_{\bar\theta,a})}$.
In particular, the two operators have the same spectral radius, and the adjoint relation exchanges their right and left eigenvalue equations.
Thus, Proposition~\ref{spectral_input_contractive}, applied with parameter $\bar\theta$, yields the claim.
\end{proof}

The effective state space is
\begin{align*}
    \Omega
    \coloneqq\begin{cases}
        [-a,1],&0<a<-\bar\theta,\quad ~(Y)
        \\(\bar\theta,1],&a\geq-\bar\theta,\quad\quad\quad (G)
    \end{cases}
\end{align*}
We briefly explain the state space. If $x\leq\bar\theta$, then $\theta x\geq1$, so no $y\in[-a,1]$ satisfies $y>\theta x$. Such a point is killed in one step, and every eigenfunction belonging to a non-zero eigenvalue must vanish there. Hence, the effective left boundary giving the two cases above is $$\ell\coloneqq\max(-a,\bar\theta).$$ 
We next show that $\lambda$ is a simple eigenvalue with a positive
eigenfunction and that all other spectral values have smaller modulus.
\begin{Lemma}\label{eigenvalue_simple_theta_smaller_minus_one_operator}
Let $\theta<-1$ and $a>0$.  $\lambda$ is algebraically simple and is the only spectral value of modulus $\lambda$ for both the killed operator restricted to $\Omega$ and the full operator $T_{\theta,a}$. Moreover, the operator restricted to $\Omega$ has an
eigenfunction $\tilde h\in L^2(\Omega,\mu_a)$ corresponding
to $\lambda$, which is strictly positive
$\mu_a$-almost everywhere on $\Omega$.
\end{Lemma}
\begin{proof}
The proof is structured as follows. We first compare the full operator with its restriction to $\Omega$ which we call $\mathcal T$.
We then compute the kernel of $\mathcal T^2$ and show that it is
strictly positive. Jentzsch's theorem gives the required properties
for $\mathcal T^2$, which we finally transfer to $\mathcal T$ and
the full operator $T_{\theta,a}$. Formally, $\mathcal T$ is defined as follows:
    \begin{align*}
        \mathcal Tf(x)\coloneqq\frac{1}{1+a}\int_{\Omega\cap\{y>\theta x\}}f(y)\d y,
        \qquad x\in\Omega.
    \end{align*}
    We regard $\mathcal T$ as an operator on $L^2(\Omega,\mu_a)$. Both $\mathcal T$ and $T_{\theta,a}$ have square-integrable kernels, since their kernels are bounded and the underlying measure spaces have finite measure. Hence, both operators are compact; see \cite[Proposition~II.4.7]{ConwayFunctionalAnalysis}.

    In case (G), so $a\geq-\bar\theta$, consider the decomposition
    \begin{align*}
        L^2([-a,1],\mu_a)
        =
        L^2([-a,\bar\theta],\mu_a)
        \oplus
        L^2(\Omega,\mu_a).
    \end{align*}
    For $f_0\in L^2([-a,\bar\theta],\mu_a)$, define
    \begin{align*}
        \mathcal T_0f_0(x)
        \coloneqq
        \frac{1}{1+a}
        \int_{[-a,\bar\theta]\cap\{y>\theta x\}}
        f_0(y)\d y,
        \qquad x\in\Omega.
    \end{align*}
    Write $f=(f_0,f_1)$, where $f_0$ is the restriction of $f$ to $[-a,\bar\theta]$ and $f_1$ is its restriction to $\Omega$. If $x\leq\bar\theta$, then $\theta x\geq1$. Since all possible values of $y$ belong to $[-a,1]$, there is no $y$ satisfying $y>\theta x$. Hence,
        $T_{\theta,a}f(x)=0$,
        $x\in[-a,\bar\theta]$.
    If $x\in\Omega$, we split the defining integral at $\bar\theta$ and obtain
    \begin{align*}
    T_{\theta,a}f(x)
    &=
    \frac{1}{1+a}
    \int_{[-a,\bar\theta]\cap\{y>\theta x\}}
    f_0(y)\d y +\frac{1}{1+a}
    \int_{\Omega\cap\{y>\theta x\}}
    f_1(y)\d y
    =(\mathcal T_0f_0)(x)+(\mathcal Tf_1)(x).
\end{align*}
Therefore,
    \begin{align*}
        T_{\theta,a}(f_0,f_1)=\bigl(0,\mathcal T_0f_0+\mathcal Tf_1\bigr).
    \end{align*}
    Thus, the full operator considered on $L^2([-a,\bar\theta],\mu_a)\oplus L^2((\bar\theta,1],\mu_a)$ has the block form
    \begin{align*}
        T_{\theta,a}
        =
        \begin{pmatrix}
            0&0\\
            \mathcal T_0&\mathcal T
        \end{pmatrix}.
    \end{align*}
    Let $\nu\neq0$ and $k\geq1$. 
    The block representation gives
    \begin{align*}
        (T_{\theta,a}-\nu I)(f_0,f_1)
        =
        \bigl(-\nu f_0, \mathcal T_0f_0+(\mathcal T-\nu I)f_1\bigr).
    \end{align*}
     Repeated application of $T_{\theta,a}-\nu I$ shows inductively that the first component of $ (T_{\theta,a}-\nu I)^k(f_0,f_1)$ is $(-\nu)^kf_0$.
    
    Consequently, if $(T_{\theta,a}-\nu I)^k(f_0,f_1)=0$,
    then $(-\nu)^kf_0=0$. Since $\nu\neq0$, this implies $f_0=0$. For vectors of the form $(0,f_1)$, the operator $\mathcal T_0$ does not contribute, and hence
    \begin{align*}
        (T_{\theta,a}-\nu I)^k(0,f_1)
        =
        \bigl(0,(\mathcal T-\nu I)^kf_1\bigr).
    \end{align*}
    It follows that
    \begin{align*}
        \ker(T_{\theta,a}-\nu I)^k
        =
        \{0\}\oplus
        \ker(\mathcal T-\nu I)^k,
        \qquad \nu\neq0,\quad k\geq1.
    \end{align*}
    Thus, the generalized eigenspaces of $T_{\theta,a}$ and $\mathcal T$ corresponding to any non-zero eigenvalue agree under the embedding $f_1\mapsto(0,f_1)$. Since both operators are compact, every non-zero spectral
    value is an eigenvalue by
    \cite[Theorem~VII.7.1, p.~214]{ConwayFunctionalAnalysis}. Therefore,
    $T_{\theta,a}$ and $\mathcal T$ have the same non-zero spectrum,
    including algebraic multiplicities.
    In case (Y), this is immediate because $\Omega=[-a,1]$ and
    $T_{\theta,a}=\mathcal T$. Together with
    \eqref{adjoint_duality_negative_theta}, this shows that $r(\mathcal T)=r(T_{\theta,a})=\lambda$.

    We now show the positivity of $\mathcal T^2$.
    For $f\in L^2(\Omega,\mu_a)$, Fubini's theorem gives, for
    $\mu_a$-almost every $x\in\Omega$,
    \begin{align*}
        \mathcal T^2f(x)
        &=
        \frac{1}{(1+a)^2}
        \int_\Omega\int_\Omega
        \mathds{1}_{\{y>\theta x\}}
        \mathds{1}_{\{z>\theta y\}}
        f(z)\d z \d y
        =\int_\Omega K_2(x,z)f(z)\d z,
    \end{align*}
    where
    \begin{align}\label{def_K_2}
        K_2(x,z)
        \coloneqq\frac{1-\max\{\ell,\theta x,\bar\theta z\}}{(1+a)^2}.
    \end{align}
    Indeed, the intermediate point $y$ must satisfy both $y>\theta x$ and $z>\theta y$, or equivalently $y>\max\{\theta x,\bar\theta z\}$.
    Since it must also belong to $\Omega$, its admissible interval is $\left(\max\{\ell,\theta x,\bar\theta z\},1\right]$. Moreover, $\theta x<1$ and $\bar\theta z<1$ for all $x,z\in\Omega$. Therefore,
    \begin{align*}
        K_2(x,z)>0,\qquad x,z\in\Omega.
    \end{align*}
    Since $r(\mathcal T)=\lambda$, by the spectral-radius formula
    \cite[Proposition~VII.3.8, p.~197]{ConwayFunctionalAnalysis},
    \begin{align*}
        r(\mathcal T^2)
        &=\lim_{n\to\infty}
        \lVert\mathcal T^{2n}\rVert^{1/n}=\left(
        \lim_{n\to\infty}
        \lVert\mathcal T^{2n}\rVert^{1/(2n)}
        \right)^2
        =r(\mathcal T)^2
        =\lambda^2.
    \end{align*}
    The measure space $(\Omega,\mu_a)$ is finite, $\mathcal T^2$ is compact, and its kernel with respect to $\mu_a$ is $(1+a)K_2$. This kernel is strictly positive $\mu_a\otimes\mu_a$-almost everywhere, so $\mathcal T^2$ is positivity improving. Hence, Jentzsch's theorem
    \cite[Theorem~V.6.6, p.~337]{SchaeferPositiveOperators} shows that $\lambda^2$ is an algebraically simple eigenvalue of $\mathcal T^2$ with a $\mu_a$-almost everywhere strictly positive eigenfunction $\tilde h\in L^2(\Omega,\mu_a)$, and every other eigenvalue has modulus strictly smaller than $\lambda^2$.
    
    If $\nu\in\sigma(\mathcal T^2)\setminus\{\lambda^2\}$, then either $\nu=0$ or, by compactness and \cite[Theorem~VII.7.1, p.~214]{ConwayFunctionalAnalysis}, $\nu$ is another eigenvalue of $\mathcal T^2$. In both cases, $|\nu|<\lambda^2$. Hence, $\lambda^2$ is the only spectral value of $\mathcal T^2$ with modulus $\lambda^2$.

    Finally, we transfer back from $\mathcal T^2$ to $\mathcal T$.
    The function $\mathcal T\tilde h$ is a non-zero eigenfunction of $\mathcal T^2$ for $\lambda^2$: indeed, $\mathcal T\tilde h=0$ would imply $\mathcal T^2\tilde h=0$, which is impossible. Since $\lambda^2$ is algebraically simple, there is a scalar $\alpha\ne 0$ such that $\mathcal T\tilde h=\alpha\tilde h$. Both $\mathcal T\tilde h$ and $\tilde h$ are $\mu_a$-almost everywhere strictly positive, so $\alpha>0$. Moreover,
    \begin{align*}
        \lambda^2\tilde h
        =\mathcal T^2\tilde h
        =\mathcal T(\alpha\tilde h)
        =\alpha^2\tilde h,
    \end{align*}
    and since $\tilde h\neq0$, this implies $\alpha^2=\lambda^2$.
    Since $\alpha>0$ and $\lambda>0$, we obtain $\alpha=\lambda$.
    The algebraic simplicity of $\lambda^2$ also implies the algebraic simplicity of $\lambda$ for $\mathcal T$. Indeed, a second eigenvector or a non-trivial Jordan chain for $\lambda$ would give the same for $\lambda^2$ under $\mathcal T^2$. More explicitly, if
    \begin{align*}
        (\mathcal T-\lambda I)g=f,\qquad
        \mathcal Tf=\lambda f,\qquad f\neq0,
    \end{align*}
    then
    \begin{align*}
    (\mathcal T^2-\lambda^2I)g
    &=(\mathcal T+\lambda I)(\mathcal T-\lambda I)g
    =(\mathcal T+\lambda I)f
    =\mathcal Tf+\lambda f
    =2\lambda f\neq0.
\end{align*}
    Thus, a non-trivial Jordan chain for $\lambda$ would give one for $\lambda^2$.
    
    It remains to exclude further spectral values of $\mathcal T$ on the circle of radius $\lambda$. Let $\nu\in\sigma(\mathcal T)$ satisfy $|\nu|=\lambda$ and $\nu\neq\lambda$. Since $\mathcal T$ is compact and $\nu\neq0$, \cite[Theorem~VII.7.1, p.~214]{ConwayFunctionalAnalysis} shows that $\nu$ is an eigenvalue. Let $f\neq0$ be a corresponding eigenvector. Then
    \begin{align*}
        \mathcal T^2f=\nu^2f.
    \end{align*}
    If $\nu=-\lambda$, then $f$ and $\tilde h$ are linearly independent eigenvectors of $\mathcal T^2$ corresponding to $\lambda^2$, contradicting its algebraic simplicity. If $\nu\neq-\lambda$, then $\nu^2\neq\lambda^2$ is another spectral value of $\mathcal T^2$ with modulus $\lambda^2$, which is also impossible. Hence, $\lambda$ is the only spectral value of $\mathcal T$ with modulus $\lambda$.
\end{proof}
We now choose a continuous representative of the eigenfunction,
normalize it, and extend it to all of $\R$.
\begin{Lemma}\label{unique_positive_eigenfunction_theta_smaller_minus_one_operator}
    Let $\theta<-1$ and $a>0$. There is a unique continuous non-negative function $h:\R\to\R$ which is strictly positive on $\Omega$, vanishes on $[-a,1]\setminus\Omega$, is normalized by $h=1$ on $[\ell\bar\theta,1]$, and satisfies \eqref{eq_eigenfunction} for every $x\in\R$.
\end{Lemma}
\begin{proof}
    By Lemma~\ref{eigenvalue_simple_theta_smaller_minus_one_operator}, there is a $\mu_a$-almost everywhere strictly positive eigenfunction $\tilde h\in L^2(\Omega,\mu_a)$ corresponding to $\lambda$. The eigenvalue equation for $\mathcal T^2$ gives, almost everywhere on $\Omega$,
    \begin{align*}
        \tilde h(x)=\frac{1}{\lambda^2}\int_\Omega K_2(x,z)\tilde h(z)\d z,
    \end{align*}
    where $K_2$ is defined in \eqref{def_K_2} in the proof of the last lemma.
    The right-hand side is continuous in $x$ and therefore defines a continuous representative of $\tilde h$. Since $K_2(x,z)>0$ and $\tilde h>0$ almost everywhere, this representative is strictly positive on $\Omega$. The identity $\mathcal T\tilde h=\lambda\tilde h$ then holds everywhere, since both sides are continuous. Thus,
    \begin{align}\label{eq_tilde_h_eigenfct}
        \lambda\tilde h(x)=\frac{1}{1+a}\int_{\max(\ell,\theta x)}^1\tilde h(y)\d y,
        \qquad x\in\Omega.
    \end{align}
    For $x\geq\ell\bar\theta$, we have $\theta x\leq\ell$, so $\tilde h$ is constant on $[\ell\bar\theta,1]$. Since this constant is strictly positive, we may rescale $\tilde h$ so that $\tilde h=1$ for $x\in[\ell\bar\theta,1]$.

    We now define $h$ on all of $\R$. Let $h=\tilde h$ on $\Omega$ and, in case (G), let $h=0$ on $[-a,\bar\theta]$. For $x\notin[-a,1]$, set
    \begin{align*}
        h(x)\coloneqq\frac{1}{\lambda(1+a)}
        \int_{\{y>\theta x\}\cap\Omega}\tilde h(y)\d y.
    \end{align*}
    We first verify the trivial regions. For $x\notin[-a,1]$, the eigenvalue equation holds by definition. In case (G), if $x\in[-a,\bar\theta]$, then $\theta x\geq1$, so both sides of the eigenvalue equation are zero. In case (Y), there are no points of $[-a,1]$ outside $\Omega$. On $\Omega$, the equation is exactly \eqref{eq_tilde_h_eigenfct}. Hence, \eqref{eq_eigenfunction} holds for every $x\in\R$. Its right-hand side is continuous in $x$, so $h$ is continuous on $\R$.

    Finally, we prove the uniqueness of $h$. Suppose that $g$ is another function with the stated
    properties. Then $g|_\Omega$ is an eigenfunction of
    $\mathcal T$ corresponding to $\lambda$. Since $\lambda$ is algebraically
    simple, $g$ is a scalar multiple of $\widetilde h$ on $\Omega$.
    The normalization on $[\ell\bar\theta,1]$ shows that this scalar is one.
    The values on $[-a,1]\setminus\Omega$ and outside $[-a,1]$ are then
    determined by the stated zero extension and the eigenvalue equation.
    Hence, $g=h$ on $\R$.
    \end{proof}
    Now, we find an explicit expression for the function $h$.
    \begin{Lemma}\label{method_of_steps_theta_smaller_minus_one_operator}
Let $\theta<-1$, $a>0$, and let $h$ be the eigenfunction from Lemma~\ref{unique_positive_eigenfunction_theta_smaller_minus_one_operator}. Recall that $\ell=\max\{-a,\bar\theta\}$. Then
\begin{align*}
    h(x)=0\quad\text{for }x\leq\bar\theta,
    \qquad
    h(x)=1\quad\text{for }x\geq\ell\bar\theta.
\end{align*}
For $\bar\theta<x<0$, we have
\begin{align}\label{neg_h_theta_neg}
    h(x)
    &=\frac{1-\theta x}{\lambda(1+a)}
    +\sum_{m=1}^{\infty}\frac{\theta^{m(2m+1)}}{[\lambda(1+a)]^{2m+1}}
    \left[
    \frac{1-\theta\ell}{(2m)!}
        (x-\ell\bar\theta^{2m})_+^{2m}-
        \frac{\theta^{2m+1}}{(2m+1)!}
        (x-\ell\bar\theta^{2m})_+^{2m+1}
    \right],
\end{align}
for $x>0$, we have
\begin{align}\label{pos_h_theta_neg}
    h(x)
    &=1+\sum_{m=1}^{\infty}\frac{\theta^{m(2m-1)}}{[\lambda(1+a)]^{2m}}
    \left[
    \frac{1-\theta\ell}{(2m-1)!}
        (\ell\bar\theta^{2m-1}-x)_+^{2m-1}+
        \frac{\theta^{2m}}{(2m)!}
        (\ell\bar\theta^{2m-1}-x)_+^{2m}
    \right],
\end{align}
and at zero,
\begin{align*}
    h(0)
    &=
    \frac12 E\left(\bar\theta, \frac{\theta\ell}{\lambda(1+a)}\right)
    +
    \frac{1-\theta\ell}{2\lambda(1+a)} E\left(\bar\theta, \frac{\ell}{\lambda(1+a)}\right).
\end{align*}
\end{Lemma}
\begin{proof}
The proof outline is the following. We first use the eigenvalue equation to determine the constant
parts of $h$ and its initial linear part. We then describe how multiplication
by $\theta$ maps the intervals $N_m$ and $P_m$ onto already
treated intervals, so that integration determines $h$
successively on $P_0,N_1,P_1,N_2,\ldots$.
Next, we check that the stated sums have the same initial
pieces and satisfy this integral equation. They therefore
agree with $h$ on every such interval. Finally, we pass to
the limit at zero and simplify the resulting series using
the characteristic equation for $\lambda$.

The constant parts for $x\leq\bar\theta$ and $x\ge\ell\bar\theta$ follow from Lemma~\ref{unique_positive_eigenfunction_theta_smaller_minus_one_operator}. 

For the remaining area, we first determine the initial linear part and the crucial differential equation. 
If $\bar\theta<x\le \ell$, the integral in the eigenvalue equation \eqref{eq_eigenfunction} only uses the constant part of $h$ since $ \ell\bar\theta\leq\theta\ell\leq1$, and therefore
\begin{align}\label{eq_h_constant_part}
    h(x)&=
    \frac{1}{\lambda(1+a)}
    \int_{\theta x}^1 h(y)\d y
    =\frac{1}{\lambda(1+a)} \int_{\theta x}^1 1\d y
    = \frac{1-\theta x}{\lambda(1+a)}.
\end{align}
For $x\in(\ell,\ell\bar\theta)$, the eigenvalue equation \eqref{eq_eigenfunction} reads as
\begin{align*}
    h(x)=\frac{1}{\lambda(1+a)}
    \int_{\theta x}^1 h(y)\d y.
\end{align*}
Since $h$ is continuous, the fundamental theorem of calculus
shows that $h$ is continuously differentiable on this interval,
with
\begin{align}\label{eq_h_dgl}
    h'(x)=\frac{-\theta}{\lambda(1+a)}h(\theta x).
\end{align}
For $m\in\N_0$, set
\begin{align*}
    N_m&\coloneqq
    [\ell\bar\theta^{2m},\ell\bar\theta^{2m+2}],\quad\text{and}\quad
    P_m\coloneqq
    [\ell\bar\theta^{2m+3},\ell\bar\theta^{2m+1}],
\end{align*}
so that it holds
$$
    \theta N_0=[\ell\bar\theta,\theta\ell];\qquad
    \theta N_m=P_{m-1},\quad m\geq1;\quad\text{and}\quad
    \theta P_m=N_m,\quad m\geq0.
$$
The intervals $N_m$ cover $[\ell,0)$, and the intervals $P_m$
cover $(0,\ell\bar\theta]$. 

Since $\theta N_0\subseteq[\ell\bar\theta,1]$, we have
$h(\theta u)=1$ for $u\in N_0$. Thus, for $x\in N_0$, using first \eqref{eq_h_dgl} and then \eqref{eq_h_constant_part},
\begin{align*}
    h(x)
    &=h(\ell)+\int_\ell^xh'(u)\d u
    =h(\ell)-\frac{\theta}{\lambda(1+a)}\int_\ell^x h(\theta u)\d u
   \\& =\frac{1-\theta\ell}{\lambda(1+a)}-\frac{\theta}{\lambda(1+a)}(x-\ell)
    =\frac{1-\theta x}{\lambda(1+a)}.
\end{align*}
On each subsequent negative interval, we integrate from its
left endpoint. On each positive interval, we integrate from
its right endpoint. This gives
\begin{align*}
    h(x)
    &=h(\ell\bar\theta^{2m})
    -\frac{\theta}{\lambda(1+a)}
    \int_{\ell\bar\theta^{2m}}^x h(\theta u)\d u,
    &&x\in N_m,\quad m\geq1,\\
    h(x)
    &=h(\ell\bar\theta^{2m+1})
    +\frac{\theta}{\lambda(1+a)}
    \int_x^{\ell\bar\theta^{2m+1}}h(\theta u)\d u,
    &&x\in P_m,\quad m\geq0.
\end{align*}
For $P_0$, the endpoint value is $h(\ell\bar\theta)=1$.
Since $\theta P_m=N_m$ and $\theta N_m=P_{m-1}$ for $m\geq1$,
the integrands are already known when the corresponding
interval is reached. The formulas therefore determine $h$
successively on $P_0,N_1,P_1,N_2,\ldots$.

We now verify that the sums in \eqref{neg_h_theta_neg} and \eqref{pos_h_theta_neg} satisfy these integral equations. On $(\ell,\ell\bar\theta)\setminus\{0\}$, the sums are locally
finite and their summands are differentiable. We may therefore
differentiate them term by term. Since $-\theta>0$, we have
    $$
    (\ell\bar\theta^{2m-1}-\theta x)_+
    =-\theta(x-\ell\bar\theta^{2m})_+\quad\text{and}\quad
    (\theta x-\ell\bar\theta^{2m})_+
    =-\theta(\ell\bar\theta^{2m+1}-x)_+.
    $$
Let $H_-$ and $H_+$ denote
the right-hand sides of \eqref{neg_h_theta_neg} and
\eqref{pos_h_theta_neg}, respectively. For $\ell<x<0$, differentiation and the first identity give
\begin{align*}
    &H_-'(x)\\
    &=
    \frac{-\theta}{\lambda(1+a)}
    +\sum_{m=1}^{\infty}
    \frac{\theta^{m(2m+1)}}{[\lambda(1+a)]^{2m+1}}
    \left[
        \frac{1-\theta\ell}{(2m-1)!}
        (x-\ell\bar\theta^{2m})_+^{2m-1}
        -\frac{\theta^{2m+1}}{(2m)!}
        (x-\ell\bar\theta^{2m})_+^{2m}
    \right]\\
    &=
    \frac{-\theta}{\lambda(1+a)}
   \! \!\left(\!
        1\!+\!\sum_{m=1}^{\infty}
        \frac{\theta^{m(2m-1)}}{[\lambda(1+a)]^{2m}}
       \! \left[
            \frac{1-\theta\ell}{(2m-1)!}
            (\ell\bar\theta^{2m-1}\!-\theta x)_+^{2m-1}
            \!+\!\frac{\theta^{2m}}{(2m)!}
            (\ell\bar\theta^{2m-1}\!-\theta x)_+^{2m}
        \right]
    \!\right)\\
    &=\frac{-\theta}{\lambda(1+a)}H_+(\theta x).
\end{align*}
For $0<x<\ell\bar\theta$, we obtain
\begin{align*}
    &H_+'(x)
    \\&=-\sum_{m=1}^{\infty}
    \frac{\theta^{m(2m-1)}}{[\lambda(1+a)]^{2m}}
    \left[
        \frac{1-\theta\ell}{(2m-2)!}
        (\ell\bar\theta^{2m-1}-x)_+^{2m-2}+\frac{\theta^{2m}}{(2m-1)!}
        (\ell\bar\theta^{2m-1}-x)_+^{2m-1} 
    \right]\\
    &=-\frac{\theta(1-\theta^2x)}{[\lambda(1+a)]^2}-\sum_{k=1}^{\infty}
    \frac{\theta^{(k+1)(2k+1)}}{[\lambda(1+a)]^{2k+2}}
    \left[
    \frac{1-\theta\ell}{(2k)!}
        (\ell\bar\theta^{2k+1}-x)_+^{2k}+\frac{\theta^{2k+2}}{(2k+1)!}
        (\ell\bar\theta^{2k+1}-x)_+^{2k+1}
    \right]\\
    &=\frac{-\theta}{\lambda(1+a)}H_-(\theta x).
\end{align*}
In the second equality, we separated the term with $m=1$,
used $1-\theta\ell+\theta^2(\ell\bar\theta-x)=1-\theta^2x$,
and did the index shift $k=m-1$ in the remaining sum. The last equality follows
from \eqref{neg_h_theta_neg}.
Thus, the two expressions in \eqref{neg_h_theta_neg} and
\eqref{pos_h_theta_neg} satisfy the same differential
equation as $h$ on $(\ell,\ell\bar\theta)\setminus\{0\}$.
All $m$ summands in \eqref{neg_h_theta_neg} vanish for
$\bar\theta<x\leq\ell\bar\theta^2$, and all summands in
\eqref{pos_h_theta_neg} vanish for $x\geq\ell\bar\theta$.
Hence,
\begin{align*}
    H_-(x)&=\frac{1-\theta x}{\lambda(1+a)},
    \quad\bar\theta<x\leq\ell\bar\theta^2,\\
    H_+(x)&=1,
    \quad  \quad \quad \quad x\geq\ell\bar\theta.
\end{align*}
Thus, both initial pieces agree with those of $h$.
The integration formulas above then give equality successively
on $P_0,N_1,P_1,N_2,\ldots$.

It remains to evaluate the common limit at zero. By Lemma~\ref{unique_positive_eigenfunction_theta_smaller_minus_one_operator}, $h$ is continuous, so it is enough to determine the limit as $x\searrow0$. Since $\ell\bar\theta^{2m-1}>0$, for every $x\geq0$ and $m\geq1$,
we have
\begin{align*}
    0&\leq(\ell\bar\theta^{2m-1}-x)_+^{2m-1}\leq(\ell\bar\theta^{2m-1})^{2m-1},\\
    0&\leq(\ell\bar\theta^{2m-1}-x)_+^{2m}\leq(\ell\bar\theta^{2m-1})^{2m}.
\end{align*}
The series obtained by substituting $x=0$ are absolutely convergent: after simplifying the powers of $\theta$, they are the even and odd parts of the series defining $E$. Dominated convergence and this even-odd decomposition give
\begin{align}\label{eq_value_h0}
    h(0)
    &=
    \frac12\left[
    E\left(\bar\theta,\frac{\theta\ell}{\lambda(1+a)}\right)+
    E\left(\bar\theta,\frac{-\theta\ell}{\lambda(1+a)}\right)
    \right]\notag\\
    &~~~+\frac{1-\theta\ell}{2\lambda(1+a)}\left[E\left(\bar\theta,\frac{\ell}{\lambda(1+a)}\right)- E\left(\bar\theta,\frac{-\ell}{\lambda(1+a)}\right)
    \right].
\end{align}
The root equation for $\lambda$ can be written, in both parameter regimes, as
\begin{align}\label{eq_root_lambda}
    E\left(\bar\theta,
    \frac{-\theta\ell}{\lambda(1+a)}\right)
    =\frac{1-\theta\ell}{\lambda(1+a)}E\left(\bar\theta,\frac{-\ell}{\lambda(1+a)}\right).
\end{align}
Indeed, if $\ell=\bar\theta$, both sides vanish by the root equation in case (G). Otherwise, $\ell=-a$, and this is precisely the root equation in case (Y). Inserting \eqref{eq_root_lambda} into \eqref{eq_value_h0} gives the stated value of $h(0)$.
\end{proof}
\begin{Remark}
The sums in both  \eqref{neg_h_theta_neg} and \eqref{pos_h_theta_neg} have only finitely many non-zero terms for every fixed $x\neq0$. Moreover, $h$ is a polynomial of degree $2m+1$ on $N_m$ and of degree $2m+2$ on $P_m$: The interval mappings show that the polynomial degrees increase by one at each step, and the leading coefficients do not vanish.
Further, for every $k\in\N$, the function $h$ is $k$ times continuously differentiable on a sufficiently small interval around zero. This follows inductively from the differential equation
\begin{align*}
    h'(x)=\frac{-\theta}{\lambda(1+a)}h(\theta x).
\end{align*}
The interval may depend on $k$. In particular, all derivatives at zero exist, and the differential equation implies
\begin{align*}
    h^{(n)}(0)
    =\left(\frac{-\theta}{\lambda(1+a)}\right)^n\theta^{n(n-1)/2}h(0).
\end{align*}
Since $h(0)>0$ and $|\theta|>1$, the Taylor series of $h$ at zero has radius of convergence zero. Thus, the polynomial pieces of $h$ cannot be combined into a convergent Taylor series around zero.
\end{Remark}
\begin{Lemma}\label{asymptotic_theta_smaller_minus_one_operator}
    Let $\theta<-1$ and $a>0$. Then there exists a constant
    $c(\theta,\lambda,a)>0$ such that
    \begin{align*}
       \lim_{n\to\infty} \frac{p_n^{a,y}(\theta)}{\lambda^{n+1}}
        = c(\theta,\lambda,a)h(y)
    \end{align*}
    uniformly in $y\in[-a,1]$.
\end{Lemma}
\begin{proof}
    Write $T\coloneqq T_{\theta,a}$. By Lemma~\ref{eigenvalue_simple_theta_smaller_minus_one_operator}, $\lambda$ is an algebraically simple eigenvalue of $T$ and the only spectral value of $T$ with modulus $\lambda$.

        Since $T$ is compact and $\lambda\neq0$,
    \cite[Theorem~VII.7.1, p.~214]{ConwayFunctionalAnalysis}
    shows that $\lambda$ is isolated in $\sigma(T)$. Hence,
    $\delta\coloneqq
        \operatorname{dist}\bigl(\lambda,
        \sigma(T)\setminus\{\lambda\}\bigr)>0$. The following construction follows
    \cite[(VII.6.9) and the discussion following it,
    p.~210]{ConwayFunctionalAnalysis}.
    Let $\Gamma\coloneqq
        \left\{z\in\C:|z-\lambda|=\frac{\delta}{2}\right\}$  oriented counterclockwise, so $\Gamma\cap\sigma(T)=\varnothing$ and its interior contains $\lambda$ but no other point of $\sigma(T)$. 
        Moreover, define
    \begin{align*}
        P\coloneqq\frac{1}{2\pi i}
        \int_\Gamma(zI-T)^{-1}\d z,
    \end{align*}
   where the integral is understood as a Bochner integral
    in $\mathcal L(\mathcal H)$ after parametrizing $\Gamma$. In Conway's notation,
    we take $A=T$ and $\Delta=\{\lambda\}$. Since $\lambda$ is isolated, we get that $\Delta$ is
    open and closed relative to $\sigma(T)$, and the operator $P$
    is precisely the Riesz idempotent $E(\Delta;A)$ defined in
    (VII.6.9), so $P^2=P$. By \cite[(VII.6.9) and the discussion following it,
    p.~210]{ConwayFunctionalAnalysis}, $P$ is a bounded projection
    which commutes with $T$. Moreover,
    $\operatorname{ran}P$ and $\ker P$ are invariant under $T$, and
    \begin{align*}
        \sigma\left(T|_{\operatorname{ran}P}\right)
        &=\{\lambda\},\qquad
        \sigma\left(T|_{\ker P}\right)
        =\sigma(T)\setminus\{\lambda\}.
    \end{align*}
    Since $Th=\lambda h$, we have
    \begin{align*}
        (zI-T)^{-1}h=\frac{1}{z-\lambda}h,
        \qquad z\in\Gamma.
    \end{align*}
    The Cauchy integral formula therefore gives
    \begin{align*}
        Ph
        =\frac{1}{2\pi i}\int_\Gamma
        \frac{1}{z-\lambda}h \d z
        =h.
    \end{align*}
    Thus, $h\in\operatorname{ran}P$.
    
    The restriction of $T$ to $\operatorname{ran}P$ is compact and has spectrum $\{\lambda\}$. Since $\lambda\neq0$, $T|_{\operatorname{ran}P}$ is invertible. Hence, the identity operator on $\operatorname{ran}P$ is compact, being the product of the bounded inverse of $T|_{\operatorname{ran}P}$ and the compact operator $T|_{\operatorname{ran}P}$. Therefore, $\operatorname{ran}P$ is finite-dimensional.

    Since $\lambda$ is algebraically simple, this space cannot contain a second
    eigenvector or a non-trivial Jordan chain corresponding to
    $\lambda$. Consequently,
    \begin{align*}
        \operatorname{ran}P=\operatorname{span}\{h\},
        \qquad TP=\lambda P.
    \end{align*}
    Set $S\coloneqq T(I-P)$. Then
    \begin{align*}
        T=\lambda P+S,
        \qquad PS=SP=0\quad\text{and}\quad T^n=(\lambda P+S)^n=\lambda^nP+S^n.
    \end{align*}
    On $\operatorname{ran}P$, the operator $S$ is zero, while on
    $\ker P$ it agrees with $T$. Hence,
    \begin{align*}
        \sigma(S)
        =\{0\}\cup\sigma\bigl(T|_{\operatorname{ker}P}\bigr)=\{0\}\cup\bigl(\sigma(T)\setminus\{\lambda\}\bigr).
    \end{align*}
    Since
    $\lambda$ is isolated in $\sigma(T)$,
    $\sigma(T)\setminus\{\lambda\}$ is compact. Moreover $\lambda$ is
    the only spectral value of $T$ with modulus $\lambda$, so $r(S)<\lambda$.
    Since $r(S)<\lambda$, the spectral-radius formula
    \cite[Proposition~VII.3.8]{ConwayFunctionalAnalysis}
    gives
    \begin{align*}
        \lim_{n\to\infty}
        \left\lVert\lambda^{-n}S^n\right\rVert^{1/n}
        =\frac{r(S)}{\lambda}<1.
    \end{align*}
    Hence, $\lambda^{-n}S^n\to0$ in operator norm and $\lambda^{-n}T^n\to P$ in operator norm.

  Since $\operatorname{ran}P=\operatorname{span}\{h\}$, we have $P\mathds{1}=Ch$ for some $C\in\C$. Moreover, $P$ is
    positive, because it is the operator-norm limit of the positive
    operators $\lambda^{-n}T^n$. Hence, $C\geq0$. Since $0\leq h\leq\lVert h\rVert_\infty\mathds{1}$ and $Ph=h$, we obtain
    \begin{align*}
       0\le \frac{Ph}{\lVert h\rVert_\infty}=\frac{h}{\lVert h\rVert_\infty}\le P\mathds{1}=Ch
    \end{align*}
    Since $h>0$ on $\Omega$, we get $C>0$, and, in $L^2([-a,1],\mu_a)$,
    \begin{align*}
        \lambda^{-n}T^n\mathds{1}\to Ch
    \end{align*}

    The preceding convergence is initially only an $L^2$-convergence and therefore does not by itself give convergence at every fixed starting point. Uniformly in $y\in[-a,1]$, the Cauchy-Schwarz inequality gives
    \begin{align*}
    |Tf(y)|
    &\leq\int_{[-a,1]\cap\{x>\theta y\}}|f(x)|\mu_a(\d x)
    \leq \lVert f\rVert_{L^2(\mu_a)}.
    \end{align*}
    and thus
    \begin{align*}
    \left|\lambda^{-n}T^n\mathds{1}(y)-Ch(y)\right|
    &=
    \frac{1}{\lambda}
    \left|
    T\left(\lambda^{-(n-1)}T^{n-1}\mathds{1}-Ch\right)(y)
    \right|
    \\&\leq
    \frac{1}{\lambda}
    \left\lVert
        \lambda^{-(n-1)}T^{n-1}\mathds{1}-Ch
    \right\rVert_{L^2(\mu_a)}
    \to0.
    \end{align*}
    Since $p_n^{a,y}(\theta)=[T^n\mathds{1}](y)$, the claim follows with $c(\theta,\lambda,a)\coloneqq\frac{C}{\lambda}$.
\end{proof}
\begin{proof}[Proof of Theorem~\ref{maintheorem_ge1} (G) and (Y)]
Lemma~\ref{unique_positive_eigenfunction_theta_smaller_minus_one_operator}
shows that $h$ is continuous and strictly positive on
$\Omega$.
Lemma~\ref{asymptotic_theta_smaller_minus_one_operator} gives
\begin{align*}
   \lim_{n\to\infty} \frac{p_n^{a,y}(\theta)}
    {c(\theta,\lambda,a)\lambda^{n+1}}
    = h(y),
    \qquad y\in[-a,1].
\end{align*}
Since $h=0$ on $[-a,1]\setminus\Omega$, all
assumptions of Lemma~\ref{Prob=Lim_int_trans_kernel} are satisfied, and the claim follows.
\end{proof}
\section{Stationary distributions}\label{sec_stationary}
We prove the existence and the representation of the stationary distribution of the limiting chain with transition kernel $p'$ claimed in Theorem~\ref{stationary_distribution}. 
\begin{proof}[Proof of Theorem~\ref{stationary_distribution}]
We first identify the left eigenfunction of $T_{\theta,a}$. If $\theta<-1$, taking adjoints in \eqref{adjoint_duality_negative_theta} gives
    $T_{\theta,a}^*=T_{\bar\theta,a}$.
Thus, the right eigenfunction of $T_{\bar\theta,a}$ given in
\cite[Theorem~1.2]{AW} is also a left eigenfunction of
$T_{\theta,a}$. We determine its formula on $\Omega$ by
checking which branches in \cite{AW} apply, and then verify that
its restriction to $\Omega$ still satisfies the left
eigenvalue equation.

First, we treat the light green region, that is, $a\geq-\bar\theta$. If $-\bar\theta\leq a\leq-1/\bar\theta$, we use case (B) of \cite[Theorem~1.2]{AW}; if $a>-1/\bar\theta$, we use case (G). For $x\in\Omega=(\bar\theta,1]$, we have
\begin{align*}
    \bar\theta\leq\bar\theta x<\bar\theta^2<1.
\end{align*}
In case (B), this lies in the range $-a\leq\bar\theta x<1$, since $-a\leq\bar\theta$. In case (G), it lies in the range $1/\bar\theta\leq\bar\theta x<1$, since $1/\bar\theta<\bar\theta$. Both cases therefore give
\begin{align*}
    g(x)=E\left(\bar\theta,\frac{-\bar\theta x}{\lambda(1+a)}\right),
    \qquad x\in\Omega.
\end{align*}

Next, we treat the light yellow region, so suppose that $0<a<-\bar\theta$, so that case (Y) of \cite[Theorem~1.2]{AW} applies. The first branch in that theorem is valid when $ -a\leq\bar\theta x\leq-\frac{a}{\bar\theta}$, or equivalently, $-\frac{a}{\bar\theta^2}\leq x\leq-\frac{a}{\bar\theta}$.
Since $-a/\bar\theta^2<-a$ and $0<-a/\bar\theta=-a\theta<1$, its intersection with $\Omega=[-a,1]$ is $[-a,-a\theta]$. On this interval, we obtain the stated deformed-exponential formula.

For $x>-a\theta$, we have $\bar\theta x<-a$, so the constant extension specified in \cite[Theorem~1.2]{AW} applies. Its value is
\begin{align*}
    g(-a\theta) = E\left(\bar\theta,\frac{a}{\lambda(1+a)}\right).
\end{align*}

These formulas inherit strict positivity on $\Omega$ from the corresponding right eigenfunctions in \cite[Theorem~1.2]{AW}. To check the left eigenvalue equation \eqref{left_eigenvalue_eq} restricted to $\Omega$, suppose that $a\geq-\bar\theta$ and $y\in\Omega$, then $\bar\theta y\geq\bar\theta\geq-a$.
Moreover, $y>\theta x$ is equivalent to $x>\bar\theta y$. Hence, no value outside $\Omega$ is part of the integral on the right-hand side of the left eigenvalue equation \eqref{left_eigenvalue_eq}. If $0<a<-\bar\theta$, then $\Omega=[-a,1]$, so no restriction is needed.

 Last, we treat the light blue region, which is $\theta>1$. Taking adjoints in \eqref{adjoint_duality_positive_theta} gives
\begin{align*}
    U_aT_{\theta,a}^*U_a^{-1}
    =T_{\bar\theta,1/a}.
\end{align*}
For $x\in\Omega=[-a,0)$, set $u=-x/a\in(0,1]$. Then $-1/a\leq\bar\theta u<1$, so we are in (B) of \cite[Theorem~1.2]{AW}. Let $h_{\bar\theta,1/a}$ be the right eigenfunction from
\cite[Theorem~1.2]{AW}. The identity above shows that
\begin{align*}
    T_{\theta,a}^*
    \bigl(U_a^{-1}h_{\bar\theta,1/a}\bigr)
    &=
    U_a^{-1}T_{\bar\theta,1/a}h_{\bar\theta,1/a}=
    \lambda U_a^{-1}h_{\bar\theta,1/a}.
\end{align*}
Therefore, the corresponding right eigenfunction of $T_{\theta,a}^\ast$ and thus the left eigenfunction of $T_{\theta,a}$ is given by
\begin{align*}
    g(x)
    &=U_a^{-1}h_{\bar\theta,1/a}(x)
    =E\left(\bar\theta, \frac{-\bar\theta(-x/a)}{\lambda(1+1/a)}\right)
    =E\left(\bar\theta, \frac{\bar\theta x}{\lambda(1+a)}\right).
\end{align*}
This function is strictly positive on $\Omega$. If $y\in[-a,0)$ and $y>\theta x$, then $x<\bar\theta y<0$. Hence, only points of $\Omega$ are part of the integral on the right-hand side of the left eigenvalue equation \eqref{left_eigenvalue_eq}.

It remains to verify stationarity. The functions $h$ and $g$ are bounded on $\Omega$ and strictly positive there. Since $\Omega$ has positive finite Lebesgue measure,
\begin{align*}
    0<\int_\Omega h(u)g(u)\d u<\infty.
\end{align*} 
For every bounded measurable function $f$, Fubini's theorem and the left eigenvalue equation \eqref{left_eigenvalue_eq} give
\begin{align*}
    \int_\Omega\int_\Omega
    f(y)p'(x,\d y)\pi(\d x)\quad
    &=
    \frac{1}
    { \int_\Omega
    h(u)g(u)\d u}
    \int_\Omega
    \frac{f(y)h(y)}{\lambda(1+a)}
    \int_{\Omega\cap\{y>\theta x\}}
    g(x)\d x \d y\\
    &=
    \frac{1}
    { \int_\Omega
    h(u)g(u)\d u}
    \int_\Omega
    f(y)h(y)g(y)\d y\\
    &=
    \int_\Omega f(y)\pi(\d y).
\end{align*}
Thus, $\pi$ is stationary. In the white regime, $h=1$, $\lambda=1$, and $p'$ is the uniform kernel, so the stationary distribution is uniform.
\end{proof}
The same argument also gives stationary distributions in the
non-expanding regime. In the following, $h_{\theta,a}$ denotes
the right eigenfunction for the indicated parameters: it is
given in \cite[Theorem~1.2]{AW} when $\theta\in[-1,1)$ and in
Theorem~\ref{maintheorem_ge1} when $|\theta|>1$.
\begin{proof}[Proof of Corollary~\ref{cor_stationary_nonexpanding}]
If $-1<a\leq0$ and $-1\leq\theta<0$, all persistence constraints
are satisfied almost surely. Thus, $h=g=1$ and
$\lambda=1$, and $\pi$ is uniform on $[-a,1]$.

    For $\theta=-1$ and $a>0$, the constraint $y>-x$ is symmetric in
    $x$ and $y$, so $T_{-1,a}^*=T_{-1,a}$. Thus, we may take
    $g=h$, and the stationary density is proportional to $h^2$.

    For $\theta=0$, the constant function $h=1$ satisfies
    \begin{align*}
        T_{0,a}h(x)
        =\frac{1}{1+a}\int_{\max\{-a,0\}}^1\d y
        =\frac{1-\max\{-a,0\}}{1+a}
        =\lambda,
        \qquad x\in[-a,1].
    \end{align*} Since $g(y)=\mathds{1}_{(0,1]}(y)$ on $[-a,1]$, we have
    \begin{align*}
        T_{0,a}^*g(y)
        &=
        \frac{\mathds{1}_{\{y>0\}}}{1+a}
        \int_{-a}^1g(x)\d x
        =\lambda g(y),
        \qquad y\in[-a,1].
    \end{align*}
    The product formula therefore gives the uniform stationary
    distribution on $[0,1]$.

    If $-1<\theta<0$ and $a> 0$, the adjoint duality \eqref{adjoint_duality_negative_theta} gives
    $T_{\theta,a}^*=T_{\bar\theta,a}$.
    Thus, $g=h_{\bar\theta,a}$ is a left eigenfunction
    corresponding to $\lambda$.

    If $0<\theta<1$ and $a>0$, the duality \eqref{adjoint_duality_positive_theta} gives
    \begin{align*}
        T_{\theta,a}^*
        =U_{1/a}T_{\bar\theta,1/a}U_{1/a}^{-1}.
    \end{align*}
    Thus, if $h_{\bar\theta,1/a}$ is the right eigenfunction
    of $T_{\bar\theta,1/a}$ from
    Theorem~\ref{maintheorem_ge1}, then
    $U_{1/a}h_{\bar\theta,1/a}$ is a left eigenfunction
    of $T_{\theta,a}$. Indeed,
    \begin{align*}
        T_{\theta,a}^*
        \bigl(U_{1/a}h_{\bar\theta,1/a}\bigr)
        &=U_{1/a}T_{\bar\theta,1/a}h_{\bar\theta,1/a}=
        \lambda U_{1/a}h_{\bar\theta,1/a}.
    \end{align*}
    Since $(U_{1/a}f)(x)=f(-x/a)$, the required function is
    \begin{align*}
        g(x)=h_{\bar\theta,1/a}(-x/a).
    \end{align*}
    Substituting the formula from
    Theorem~\ref{maintheorem_ge1} (B) gives the stated
    expression for $g$.

Last, if $0<\theta<1$ and $-1<a<0$, we distinguish the following two cases.
If $\theta\leq-a$, the stated sum gives $g=1$
and $\lambda=1$. Otherwise, termwise differentiation
and the same index shift as above give
\begin{align*}
    g'(y)=\frac{g(y/\theta)}{\theta\lambda(1+a)},
    \qquad -a<y<\theta,
\end{align*}
while $g=1$ on $[\theta,1]$. Consequently,
\begin{align*}
    F(y)\coloneqq
    \lambda g(y)
    -\frac{1}{1+a}
    \int_{-a}^{\min\{y/\theta,1\}}g(x)\d x
    =\lambda g(y)
    -T^*_{\theta,a}g(y)
\end{align*}
is constant on $[-a,1]$. Using $T_{\theta,a}h=\lambda h$
and Fubini's theorem, we obtain
\begin{align*}
    \int_{-a}^1h(y)F(y)\d y
    =\lambda\int_{-a}^1hg
    -\int_{-a}^1g(x)T_{\theta,a}h(x)\d x
    =0.
\end{align*}
Since $h>0$, the constant function $F$ is zero, proving the left
eigenvalue equation.

To verify non-negativity of $g$, the triangle inequality gives
$T_{\theta,a}^*|g|-\lambda|g|\geq |T_{\theta,a}^*g|-\lambda|g|=0$.
By applying Fubini and $T_{\theta,a}h=\lambda h$ as before, the integral against $h$ is zero:
\begin{align*}
    \int_{-a}^1h(y)
    \bigl(T_{\theta,a}^*|g|(y)-\lambda|g(y)|\bigr)\d y
    &=\int_{-a}^1|g(x)|
    \bigl(T_{\theta,a}h(x)-\lambda h(x)\bigr)\d x=0.
\end{align*}
Since $h>0$ and the difference is non-negative and continuous, it vanishes
everywhere, so $T_{\theta,a}^*|g|=\lambda|g|$. Evaluating both eigenvalue equations at one,
where $g(1)=1$, gives
\begin{align*}
    \int_{-a}^1|g(x)|\d x
    =\lambda(1+a)
    =\int_{-a}^1g(x)\d x.
\end{align*}
Thus, $g\geq0$.

    In all cases, $h$ and $g$ are bounded and non-negative,
and their product is strictly positive on an interval
of positive length. Therefore,
\begin{align*}
    0<\int_{-a}^1h(u)g(u)\d u<\infty.
\end{align*}
The same Fubini calculation as in the proof of
Theorem~\ref{stationary_distribution}, now integrating
over $[-a,1]$, proves stationarity.
\end{proof}

\section*{\small Declaration of AI use}
The authors formulated the
mathematical questions and developed the mathematical ideas, results and their proofs. ChatGPT was used for additional checks
of the proofs, to help locate two references on operator theory
\cite{ConwayFunctionalAnalysis,SchaeferPositiveOperators}, improve the wording, and suggest TikZ commands for minor
adjustments to Figure~\ref{fig:map}.

\end{document}